\documentclass[onefignum,onetabnum,final]{siamart220329}

\usepackage{amsmath,amssymb,enumitem}
\usepackage{stmaryrd}
\usepackage{txfonts}

\newtheorem{rmk}{Remark}

\newcommand{\vi}{\varv}
\newcommand{\R}{\mathbb{R}}
\newcommand{\N}{\mathbb{N}}

\newcommand{\lapv}{\Delta_{\vi}}
\newcommand{\Dv}{D_{\vi}}
\newcommand{\vDx}{\vi \cdot D_{x}}
\newcommand{\divv}{\text{div}_{\vi}}
\newcommand{\pt}{\partial_t}
\newcommand{\Lt}{L^{2}([0,T]\times \R^d\times \R^d)}
\newcommand{\intr}{\int_{\R^d\times \R^d}}
\newcommand{\me}{m^{\epsilon}}
\newcommand{\ue}{u^{\epsilon}}
\newcommand{\gep}{g^{\epsilon}}

\newcommand{\un}{u^n}

\newcommand{\phil}{\phi_{\lambda}}

\newcommand{\Phia}{\Phi_{\alpha}}
\newcommand{\Phiak}{\Phi_{\alpha,k}}

\title{Weak and renormalized solutions to a hypoelliptic Mean Field Games system.}
\author{Nikiforos Mimikos-Stamatopoulos\thanks{University of Chicago, Department of Mathematics, Chicago, IL}}

\makeatletter
\def\namedlabel#1#2{\begingroup
    #2%
    \def\@currentlabel{#2}%
    \phantomsection\label{#1}\endgroup
}
\makeatother
\begin{document}
\maketitle
\begin{abstract}
We study the well-posedness of a degenerate, hypoelliptic Mean Field Games system with local coupling and Hamiltonians which are either Lipschitz or grow quadratically in the gradient. In the former case, we prove the existence and uniqueness of weak solutions while in the latter we study the same question for renormalized solutions. Our approach relies on the kinetic regularity of hypoelliptic equations obtained by Bouchut and the work of Porretta on the existence and uniqueness of renormalized solutions for the non-degenerate Mean Field Game system.
\end{abstract}
\begin{keywords}
MFG, Hypoelliptic
\end{keywords}
\begin{AMS}
 35Q91 (Primary) 91A13 35H10 (Secondary)
\end{AMS}
\section*{Introduction}
We establish the well posedness (existence and uniqueness) of solutions of the local, hypoelliptic Mean Field Games system (MFG for short)
\begin{equation}
\label{MFGequation}
\begin{cases}
-\pt u-\lapv u+\vDx u+H(\Dv u)=F(t,x,\vi,m(t,x,\vi))\text{ in }(0,T)\times \R^d\times \R^d,\\
\pt m-\lapv m-\vDx m-\divv(m H_p(\Dv u))=0, \text{ in }(0,T)\times \R^d\times \R^d,\\
u(T,x,\vi)=G(x,\vi,m(T,x,\vi)), m(0,x,\vi)= m_0(x,\vi).
\end{cases}
\end{equation}
The Hamiltonian $H:\R^d\rightarrow \R$ is convex, the coupling term $F:[0,T]\times \R^d\times \R^d\times\R\rightarrow \R$ as well as the terminal cost function $G:\R^d\times \R^d\times \R\rightarrow \R$ are increasing in $m$, and $m_0:\R^d\times \R^d\rightarrow \R$ is a given probability density.\par   
Systems like \eqref{MFGequation} formally describe the equilibrium of an $N$-player game, when $N$ tends to infinity, of indistinguishable players, where each player makes decisions based on the distribution of the other co-players. In this setup, it is natural to interpret $x\in \R^d$ as the position and $\vi\in \R^d$ as the velocity of such players. More precisely, the players control their accelaration in order to minimize the cost introduced by the coupling $F$ and the Hamiltonian $H$, which leads to the Hamilton-Jacobi-Bellman equation (HJB for short). The optimal feedback is then given by the vector field $-(\vi,D_pH(\Dv u))$, under which, their distribution changes according to the degenerate Fokker-Planck equation (FP for short). As far as applications are concerned, we refer to the flocking model in Carmona and Delarue \cite{carmona2018probabilistic}, and for a first order system we refer to Bardi and Cardaliaguet \cite{bardi2021convergence}, Griffin and Meszaros \cite{griffin2022variational} and Achdou, Mannucci, Marchi and Tchou \cite{achdou2020deterministic}. Finally, we mention that the general form of \eqref{MFGequation} is reminiscent of Boltzmann-type equations, which have been investigated in the MFG context by Burger, Lorz, Wolfram \cite{burger2016balanced} in a setting different to the one used in this paper.\par  
MFG were introduced by Lasry and Lions in \cite{lasry2006jeux}, \cite{lasry2006_2jeux}, \cite{lasry2007mean} and, in a special case, by Huang, Caines, Malhame \cite{huang2006large}. Although there has been extensive study of non-degenerate second-order mean field games, with a local or non-local coupling, less has been done in the degenerate setting, an example of the latter being hypoelliptic MFG. In this setting, when the degeneracy is a sum of squares, Dragoni and Feleqi studied in \cite{dragoni2018ergodic} the ergodic problem; see also Feleqi, Gomes and Tada \cite{feleqi2020hypoelliptic}. When $H(p)=\frac{1}{2}|p|^2$, Camilli in \cite{camilli2021quadratic}, obtained, using the Hopf-Cole transformation, weak solutions to \eqref{MFGequation} with uncoupled terminal data. We remark that the assumptions of Camilli appear almost complementary to the ones in this paper, as the existence of solutions in \cite{camilli2021quadratic} is established for terminal data that have to be unbounded since they need to be superquadratic. For results in the case of non-H\"ormander degenerate systems, we refer to Cardaliaguet, Graber, Porretta and Tonon in \cite{cardaliaguet2015second}, who study, using a variational approach, degenerate MFG systems, for Hamiltonians with super-linear growth and no coupling on the terminal data of the HJB equation.\par
Our goal is to show existence and uniqueness for quadratic and Lipschitz Hamiltonians, under similar assumptions as that of Porretta in \cite{porretta2015weak}, where existence and uniqueness of renormalized solutions was established in the non-degenerate setting. We work with two different types of Hamiltonian $H$, that is, with linear or quadratic growth. Furthermore, the degeneracy is not a sum of squares, that is, $L$ is not of the form $L:=\sum\limits_{i,j}^k a_{ij}X_iX_j$, for some vector fields $X_i$ satisfying H\"ormanders condition . In the context of hypoelliptic operators, the degenerate operator $L:=\pt -\lapv +\vDx$ is the simplest and historically the first one to be studied.\\
The first result addresses the case of a Lipschitz Hamiltonian, whereas the latter the case of quadratic Hamiltonian. 
\begin{theorem}\label{ExistenceTheorem}
Assume that $H,F,G,$ and $m_0$ satisfy \ref{AssumptionHamilt},\ref{AssumptionCoupling},\ref{AssumptionTerminal}, and \ref{AssumptionInitial}. Then, there exists a unique weak solution $(u,m)$ of \eqref{MFGequation}, according to Definition \ref{WeakSolutionsDef}. Moreover, there exists a constant $C>0$, such that,
\[\|-\pt u+\vDx u\|_{L^2([0,t]\times \R^d\times \R^d)}+\|\lapv u\|_{L^2([0,t]\times \R^d\times \R^d)}\]
\[+\|-\pt m+\vDx m\|_{L^2([0,t]\times \R^d\times \R^d)}+\|\lapv m\|_{L^2([0,t]\times \R^d\times \R^d)}\leq \frac{C}{T-t}.\]
Furthermore, if $F$ also satisfies \ref{AssumptionRegularity}, there exists a constant $C=C(F,G,H,T,m_0)>0$, such that
\[\sup\limits_{t\in [0,T]}\|m(t)\|_2+\sup\limits_{t\in [0,T]}\|Dm(t)\|_2+\|D_{\vi,\vi}^2 m\|_2+\|\Dv D_x m\|_2\leq C,\]
and
\[\sup\limits_{t\in [0,T]}\|u(t)\|_2+\sup\limits_{t\in [0,T]}\|Du(t)\|_2+\|D_{\vi,\vi}^2 u\|_2+\|\Dv D_x u\|_2\leq C.\]
\end{theorem}
The second result is about renormalized solutions as in Definition \ref{DefinitionOfRenormalizedToMFGSystem}.
\begin{theorem}\label{QuadraticExistenceTheorem}
Assume that $H,F,G,$ and $m_0$ satisfy \ref{Assumption2onH},\ref{F_LAssum},\ref{G_LAssum}, and \ref{AssumptionInitial}. Then, there exists a unique pair $(u,m)$, of renormalized solutions of the MFG system \eqref{MFGequation}. Furthermore, assume that $F,G$ are only functions of $m$. Then, there exists a constant $C=C(m_0,F,G,T)>0$, such that 
\[\intr G'(m(T,x,\vi))|Dm(T,x,\vi)|^2dxd\vi +\int_0^T\intr F'(m(t,x,\vi))|Dm(t,x,\vi)|^2\]
\[+m\sum\limits_{k=1}^{2d}m\Dv u_{k}H_{pp}(\Dv u)\Dv u_k dxd\vi \leq C.\]
\end{theorem}
The existence of a solution, in the case of Lipschitz Hamiltonians, is established using a Schauder fixed point theorem as follows. Fix a probability density $m_0$. Given $\mu\in X:=C([0,T];L^2(\R^d\times \R^d)),$ let $u^{\mu}\in C([0,T];L^2(\R^d\times \R^d))$, with $\Dv u\in L^2([0,T]\times \R^d\times \R^d),$ be the unique, distributional solution of 
\begin{equation*}
\begin{cases}
-\pt u-\lapv u+\vDx u+H(\Dv u)=F(t,x,\vi,\mu) \text{ in }(0,T)\times \R^d\times \R^d,\\
u(T,x,\vi)=G(\mu(T,x,\vi)) \text{ in }\R^d\times \R^d,
\end{cases}
\end{equation*}
and $m$ the unique distributional solution of 
\begin{equation*}
    \begin{cases}
    \pt m-\lapv m-\vDx m-\divv(mD_pH(\Dv u^{\mu}))=0 \text{ in }(0,T)\times \R^d\times \R^d,\\
    m(0,x,\vi)=m_0(x,\vi) \text{ in }\R^d\times \R^d.
    \end{cases}
\end{equation*}
Set $\Phi(\mu)=m$. We need to show that $\Phi$ is $X-$valued, continuous, and compact. The two aforementioned properties follow easily once we show that $\Phi(m)\in L^{\infty}$ with appropriate bounds. Compactness does not follow immediately, because of the degenerate $x-$direction. To work with that, we localize in time the results in Bouchut \cite{bouchut2002hypoelliptic}.\par
 For Theorem \ref{QuadraticExistenceTheorem}, we rely on the work in \cite{porretta2015weak} and mostly adapt the arguments in the hypoelliptic setting. In particular, given a Hamiltonian $H$ with quadratic growth (exact assumptions are given later), we consider a sequence of Lipschitz pointwise-approximations and the corresponding solutions provided by Theorem \ref{ExistenceTheorem} and show compactness in the appropriate spaces. The main technical difficulties and deviations from \cite{porretta2015weak} are the gradient estimates in hypoelliptic equations with $L^1-$data, which are briefly described next. Let $H^{\epsilon}$ be a suitable pointwise Lipschitz approximation of a quadratic Hamiltonian $H$ and $(\me,\ue)$ the corresponding weak solutions. In order to show that there exists a limit which is a renormalized solution, we must show the convergence (up to a subsequence) of $\ue,\me$ in $L^1([0,T]\times \R^d\times \R^d)$ and of the gradients $\Dv(\ue \wedge k) ,\Dv (\me\wedge k )$ of the truncations in $L^2([0,T]\times \R^d\times \R^d)$. The compactness of $\ue$ in $L^1$ follows by the results of DiPerna and Lions in \cite{diperna1988fokker}, while the convergence of the gradients is due to an appropriate transformation similar to the one used by Porretta in \cite{porretta1999existence} and the references therein. This important transformation is studied in the Appendix. Finally, for the FP equation, the crucial bound as pointed out in \cite{porretta2015weak} is that, for some independent of $\epsilon$, $C>0$, 
\begin{equation}
\label{PorrettaEstimate}
\|\me |H_p^{\epsilon}(\Dv \ue)|^2\|_1\leq C. 
\end{equation} 
 This estimate is crucial in the following way. If $\me$ is a solution to the FP equation \eqref{MFGequation}, a priori, the best independent of $\epsilon$ estimate for $\me H_p^{\epsilon}(\Dv \ue)$ is in $L^1([0,T]\times \R^d\times \R^d)$. However, to obtain fractional gradient estimates we need bounds in $L^{r}$ for some $r>1$. The main observation that allows us to obtain this under condition \eqref{PorrettaEstimate}, is the following: Due to hypoellipticity, higher integrability of $\me H_p^{\epsilon}(\Dv \ue)$ should yield higher integrability for $\me$, while under condition \eqref{PorrettaEstimate} higher integrability of $\me$ should also yield higher integrability for $\me H_p^{\epsilon}(\Dv \ue)$. We show that it is possible to combine the above gains and obtain higher integrability with bounds independent of $\epsilon$ and therefore use the results from \cite{bouchut2002hypoelliptic}.
\subsection{Organization of the Paper}
In section 1, we state all the assumptions and definitions used throughout the paper. In section 2, we study the backwards HJB and FP equations with $L^2-$terminal/initial data respectively. The main estimates come from Theorems \ref{Thm1.2FiniteTime} and \ref{Thm1.3FiniteTime}. We also obtain results regarding the hypoelliptic FP equation and, in particular, we establish fractional gradient bounds. Finally, we establish Theorem \ref{ExistenceTheorem}. Section 3 is devoted to the proof of Theorem \ref{QuadraticExistenceTheorem}. Finally, in the appendix (section 4) we show an important technical result for the hypoelliptic HJB equation and we give the statements of the theorems we will use from \cite{bouchut2002hypoelliptic}.
\subsection{Notation and Terminology}
Throughout the paper, $d\in \N:=\{1,\cdots,\infty\}$, $T>0$ is the terminal time, $t\in [0,T]$ is the time variable, $x\in \R^d$ and $\vi,v\in \R^d$, and vectors in $[0,T]\times \R^d\times \R^d$ always appear in the order $(t,x,\vi)$. For $p\in [1,\infty]$, $L^p([0,T]\times \R^d\times \R^d)_+$ and $L^p(\R^d\times \R^d)_+$, are the non-negative functions of $L^p([0,T]\times \R^d\times \R^d)$ and $L^p(\R^d\times \R^d)$ respectively. For $s>0, W^{s,p}(\R^d\times \R^d)$ is the usual fractional Sobolev space and $D^s=(-\Delta )^{s/2}$, we refer for example to \cite{demengel2012fractional} for the definition of fractional Sobolev spaces. If $\phi=\phi(t,x,\vi):[0,T]\times \R^d\times \R^d\rightarrow \R$ or $\phi=\phi(x,\vi):\R^d\times \R^d\rightarrow \R$, we write $D^2\phi=D_{x,\vi}^2\phi$, for the hessian in the space variables, $\lapv \phi:=\sum\limits_{i=1}^d \partial_{\vi_i \vi_i}\phi$, $D\phi :=(D_x\phi,D_{\vi} \phi)$ and $\divv(\phi):=\sum\limits_{i=1}^d \partial_{\vi_i}\phi$. For a function $F(t,x,\vi,m):[0,T]\times \R^d\times \R^d\times \R\rightarrow \R$ or $G(x,\vi,m):\R^d\times \R^d\times \R\rightarrow \R$, we use the notations $D_{(x,\vi)}F=(\partial_{x_1}F,\cdots,\partial_{x_d}F,\partial_{\vi_1}F,\cdots,\partial_{\vi_d}F),\, F_m=\partial_m F$, and similarly for $G$. Throughout the paper when we reference a standard sequence of mollifiers $\rho_n:\R^d\times \R^d\rightarrow [0,\infty)$ we mean that $\rho_n(x,\vi):=n^{2d}\rho(\frac{x}{n},\frac{\vi}{n})$ where $\rho\in C_c^{\infty}(\R^d\times \R^d)$, such that $\rho\geq 0$ and $\intr \rho(x,\vi)dxd\vi =1$. Moreover in all the proofs constants are subject to change from line to line and they only depend on the quantities/functions stated in the statement of the result. Finally, we will often use the terminology dimensional constant referring to a constant that only depends on the dimension.
\section{Assumptions/Definitions}
We split this section in two subsections, one for Lipschitz Hamiltonian and one for quadratic.
\subsubsection{Lipschitz Hamiltonian and weak solutions}
As far as the data are concerned, we assume the following, for the case of Lipschitz Hamiltonian:
\begin{itemize}
    \item[\namedlabel{AssumptionHamilt}{[H1]}](Lipschitz Hamiltonian) The Hamiltonian $H:\R^d\rightarrow \R,$ is $C^1(\R^d),$ convex, $H\geq 0$, $H(0)=0$, and there exists an $L_H>0$, such that,
    \begin{equation}\tag{H1.1}
        |H(p_2)-H(p_1)|\leq L_H|p_2-p_1|\,\, \text{  for all }p_1,p_2\in \R^d.
    \end{equation}
    \item[\namedlabel{AssumptionCoupling}{[F1]}] (Coupling term) The coupling term $F=F(t,x,\vi,m):[0,T]\times \R^d\times \R^d\times\R\rightarrow \R,$ is continuous, strictly increasing and locally Lipschitz in $m$, that is, for all $L>0$, there exists a constant $c_L>0$ such that $|F(t,x,\vi,m_2)-F(t,x,\vi,m_1)|\leq c_L|m_2-m_1|$ for all $0\leq m_1,m_2\leq L$, and $F(t,x,\vi,0)\in L^2([0,T]\times \R^d\times \R^d)$. Finally, we assume that $F\geq 0$.
    \item[\namedlabel{AssumptionTerminal}{[G1]}](Terminal data for $u$) The coupling term $G=G(x,\vi,m): \R^d\times \R^d\times\R\rightarrow \R,$ is continuous, strictly increasing and locally Lipschitz in $m$ (in the same sense as $F$ above), and $G(x,\vi,0)\in L^2(\R^d\times \R^d)$. Finally, we assume that $G\geq 0$.
    \item[\namedlabel{AssumptionInitial}{[M1]}](Initial density) The initial density $m_0:\R^d\times \R^d\rightarrow \R $, satisfies $m_0\in L^{\infty}([0,T]\times \R^d\times \R^d)_+, \sqrt{m_0}\in L^1(\R^d\times \R^d), (|x|^2+|\vi|^4)m_0\in L^1(\R^d\times \R^d)$, $\log(m_0)\in L_{loc}^1(\R^d\times \R^d)$, $Dm_0\in L^{2}(\R^d\times \R^d)$ and $\int_{\R^d\times \R^d}m_0(x,\vi)dxd\vi =1$.
    \item[\namedlabel{AssumptionRegularity}{[R1]}](Regularity) Assume that $F,G$ satisfy \ref{AssumptionCoupling},\ref{AssumptionTerminal} and that for every $L>0$, there exists a $c_0=c_0(L)>0,$ such that,
\[c_0\leq |F_m(t,x,\vi,m)|,|G_m(x,\vi,m)|, \text{ for all }(t,x,\vi,m)\in [0,T]\times \R^d\times \R^d\times [0,L].\]
Furthermore, we assume that there exists a constant $C>0$, such that,
\[|D_{(x,\vi)}F(t,x,\vi,m)|+|D_{(x,\vi)}G(x,\vi,m)|\leq C|m| \text{ for all }(t,x,\vi,m)\in [0,T]\times \R^d\times \R^d\times \R.\]
\end{itemize}
\begin{rmk}
    We note that assumption \ref{AssumptionInitial} implies in particular that $m_0\log(m_0)\in L^1(\R^d\times \R^d)$.
\end{rmk}
Next we state the definition of a weak solution.
\begin{definition}\label{WeakSolutionsDef}
Assume that $H,G,F,$ and $m_0$ satisfy \ref{AssumptionHamilt},\ref{AssumptionCoupling},\ref{AssumptionTerminal} and \ref{AssumptionInitial}. A pair $(u,m)\in \Lt\times \Lt$ is a weak solution of the system \eqref{MFGequation}, if 
\[u\in C([0,T];L^2(\R^d\times \R^d), \Dv u\in \Lt),\]
\[m\in C([0,T];L^2(\R^d\times \R^d)), \Dv m\in L^2, D_x^{1/3}m\in \Lt, m\in L^{\infty}([0,T]\times \R^d\times \R^d),\]
the system \eqref{MFGequation} holds in a distributional sense.

\end{definition}
\subsubsection{Quadratic Hamiltonian and renormalized solutions}
For the case of a quadratic Hamiltonian $H$, we assume the following:
\begin{itemize}
    \item[\namedlabel{Assumption2onH}{[H2]}] For the Hamiltonian $H:\R^d\rightarrow \R$ we assume that it is convex, continuous and there exist constants $c>0,C>0$ such that, for all $p\in \R^d$,
\begin{equation}\tag{H2.1}
\label{Hamilt1}
0\leq H(p)\leq C|p|^2,
\end{equation}
\begin{equation}\tag{H2.2}
\label{Hamilt2}
H_p(p)\cdot p-H(p)\geq cH(p),
\end{equation}
\begin{equation}\tag{H2.3}
\label{Hamilt3}
|H_p(p)|\leq C|p|.
\end{equation}
\item[\namedlabel{F_LAssum}{[F2]}]
For the coupling term $F:[0,T]\times \R^d\times \R^d\times \R\rightarrow \R$, we assume that it satisfies \ref{AssumptionCoupling} and with bounds that possibly depend on $L>0$, one of the following hold:
\begin{equation}\tag{F2.1}\label{PorrettaAssumptionF}
    f_L(t,x,\vi):=\sup\limits_{m\in [0,L]}F(t,x,\vi,m)\in L^1(\R\times \R^d\times \R^d),
\end{equation}
\begin{equation}\tag{F2.2}\label{F/mAssumption}
    f_L(t,x,\vi):=\sup\limits_{m\in [0,L]} F(t,x,\vi,m)/m\in L^{\infty}(\R\times \R^d\times \R^d).
\end{equation}
\item[\namedlabel{G_LAssum}{[G2]}]
For the coupling term $G:\R^d\times \R^d\times \R\rightarrow \R$, we assume that it satisfies \ref{AssumptionTerminal} and with bounds that possibly depend on $L>0$, one of the following hold:
\begin{equation}\tag{G2.1}\label{PorrettaAssumptionG}
    g_L(x,\vi):=\sup\limits_{m\in [0,L]}G(x,\vi,m)\in L^1( \R^d\times \R^d),
\end{equation}
\begin{equation}\tag{G2.2}\label{G/mAssumption}
    g_L(x,\vi):=\sup\limits_{m\in [0,L]}G(x,\vi,m)/m\in L^{\infty}( \R^d\times \R^d).
\end{equation}
\end{itemize}
\begin{rmk}\label{RemarkUniformIntegralibity}
The above conditions on $F,G$ yield that if \eqref{PorrettaAssumptionF} and \eqref{PorrettaAssumptionG} hold, then
\begin{equation*}
F(x,\vi,t,m)\leq f_L(t,x,\vi)+\frac{m}{L}F(t,x,\vi,m),\,\,G(x,\vi,m)\leq g_L(x,\vi)+\frac{m}{L}G(x,\vi,m),
\end{equation*}
for every $m\geq 0,L>0$. While if \eqref{F/mAssumption} and \eqref{G/mAssumption} hold, then,
\[ F(t,x,\vi,m)\leq f_L(t,x,\vi)m+\frac{m}{L}F(x,\vi,m),G(x,\vi,m)\leq g_L(x,\vi)m+\frac{m}{L}G(x,\vi,m).\]
Conditions \eqref{PorrettaAssumptionF},  \eqref{PorrettaAssumptionG} do not allow for $F,G$ to depend only on $m$ due to the unbounded domain, while conditions \eqref{F/mAssumption}, \eqref{G/mAssumption} do allow for dependence only on $m$. Typical examples for the coupling are of the form
\begin{equation*}
    F(t,x,\vi,m)=a(t,x,\vi)h_1(m)+h_2(m)
\end{equation*}
where for assumption \eqref{AssumptionCoupling} we need $h_2(0)=0$, $h_1\geq 0$, strictly increasing and locally Lipschitz continuous and finally $a\geq 0$, $a\in L^2\cap L^{\infty}$ and continuous. For assumption \eqref{F_LAssum} we need to also assume that
\begin{itemize}
    \item In the case of \eqref{PorrettaAssumptionF}, $a\in L^1$ and $h_2(m)=0$.
    \item While for the case of \eqref{F/mAssumption}, we may also impose $a\in L^{\infty}$ and that $h_1(m)=m^{q_1}, h_2(m)=m^{q_2}$ for some $q_1,q_2\in [1,\infty)$.
\end{itemize}
\end{rmk}
Next, we define renormalized solutions for equations of the form
\begin{equation}
\label{RenormDefnFP}
\begin{cases}
\pt m-\lapv m-\vDx m-\divv(mb)=0 \text{ in }(0,T)\times \R^d\times \R^d,\\
m(0,x,\vi)=m_0(x,\vi) \text{ in }\R^d\times \R^d,
\end{cases}
\end{equation}
where $b:[0,T]\times \R^d\times \R^d\rightarrow \R$, $m_0:\R^d\times \R^d$, and equations of the form
\begin{equation}
\label{RenormHJB}
\begin{cases}
-\pt u-\lapv u+\vDx u+H(\Dv u)=f \text{ in }(0,T)\times \R^d\times \R^d,\\
u(T,x,\vi)=g(x,\vi) \text{ in }\R^d\times \R^d.
\end{cases}
\end{equation}
\begin{rmk}
Regarding our notation, in the rest of the paper, we will follow the convention that capital letters $F,G$ are be used when referring to the MFG system, while lower case letters $f,g$ will be used for general HJB equations.
\end{rmk}
Our definitions are in the same spirit as in \cite{porretta2015weak}.
\begin{definition}\label{RenormFPDefinition}
Let $m\in C([0,T];L^1(\R^d\times \R^d)_+)$ and $b:[0,T]\times \R^d\times \R^d\rightarrow \R$, such that $m|b|^2\in L^1([0,T]\times \R^d\times \R^d)$. We say that $m$ is a renormalized solution of \eqref{RenormDefnFP}, if
\[\lim\limits_{n\rightarrow \infty}\frac{1}{n}\int_{n<m<2n}|\Dv m|^2dxd\vi dt=0,\]
and for each $S\in W^{2,\infty}(\R), S(0)=0,$ the function $S(m)$ satisfies in the distributional sense,
\[\pt S(m)-\lapv S(m)-\vDx S(m)-\divv(S'(m)mb)+S''(m)|\Dv m|^2+S''(m)mb\Dv m=0,\]
\[S(m)(0)=S(m_0).\]

\end{definition}
\begin{definition}\label{RenormHJBDefinition}
Let $u\in C([0,T];L^1(\R^d\times \R^d)_+)$, with $\Dv u\in L^2([0,T]\times \R^d\times \R^d)$, $f\in L^1([0,T]\times \R^d\times \R^d), g\in L^1(\R^d\times \R^d)$. We say that $u$ is a renormalized solution of \eqref{RenormHJB}, if
\[\lim\limits_{n\rightarrow \infty}\frac{1}{n}\int_{n<m<2n}|\Dv u|^2dxd\vi dt=0,\]
and for each $S\in W^{2,\infty}(\R^d), S(0)=0,$ the function $S(u)$ satisfies  in the distributional sense,
\begin{equation*}
    -\pt S(u)-\lapv S(u)+\vDx S(u)+S'(u)H(\Dv u)=S'(u)f,\,\, S(u(T))=S(g).
\end{equation*}
\end{definition}
\begin{definition}\label{DefinitionOfRenormalizedToMFGSystem}
Assume that $H,G,F,$ and $m_0$ satisfy \ref{Assumption2onH},\ref{G_LAssum},\ref{F_LAssum}, and \ref{AssumptionInitial}. A pair $(m,u)\in C([0,T];L^1(\R^d\times \R^d)_+)\times C([0,T];L^1(\R^d\times \R^d)_+)$, is a renormalized solution of the MFG system \eqref{MFGequation}, if $m,u$ are renormalized solutions to the corresponding equations according to Definitions \eqref{RenormDefnFP}, \eqref{RenormHJB}, respectively.
\end{definition}
\begin{rmk}
In general the notions of renormalized and distributional solutions are distinct. However under suitable conditions we may show they are equivalent. We do not explore this direction in the present work, although it should follow with similar methods as in the non degenerate case, see Porretta \cite{porretta2015weak} and for results on the whole space Porretta \cite{porretta2017weak}.
\end{rmk}
\section{The well possedness in the case of Lipschitz Hamiltonian}
All the equations in the rest of the section should be understood in the distributional sense, unless stated otherwise. We divide this section in four parts. In the first two we study the HJB equation and the FP equation separately, in the third section we use these bounds to obtain weak solutions to the MFG problem, and in the last part we show a regularity result for these weak solutions.
\subsection{Estimates for the Hamilton-Jacobi-Bellman equation}
\begin{theorem}\label{EstimatesonHJB}
Let $g\in L^2(\R^d\times \R^d)\cap L^{\infty}( \R^d\times \R^d)_+$, $f\in C([0,T];L^2(\R^d\times \R^d))\cap L^{\infty}([0,T]\times \R^d\times \R^d)_+$, and a Hamiltonian $H:\R^d\rightarrow \R$, which satisfies \ref{AssumptionHamilt}. Then, there exists a unique solution $u\in C([0,T];L^2(\R^d\times \R^d)),$ with $\Dv u\in L^2([0,T]\times \R^d\times \R^d)$ of \eqref{RenormHJB}. Furthermore, there exists a $C=C(T,\text{Lip}_H)>0$, such that
\[\sup\limits_{t\in [0,T]}\|u(t)\|_2+\|\Dv u\|_2\leq C(\|g\|_2+\|f\|_2)\]
and for each $t\in [0,T]$,
\[\|\pt u-\vDx u\|_{L^2([0,t]\times \R^d\times \R^d)}+\|\lapv u\|_{L^2([0,t]\times \R^d\times \R^d)}\leq \frac{C}{T-t}\Big(\|f\|_2+\|g\|_2 \Big).\]
Finally, there exists a constant $C=C(T,d,\|f\|_{\infty},\|g\|_{\infty})>0$, such that $\|u\|_{\infty}\leq C,$
    in particular $C$ does not depend on the Lipschitz constant of the Hamiltonian $H$.
\end{theorem}
\begin{proof}
First we address the issue of existence. Consider the Banach space $X:=\{v\in L^2([0,T]\times \R^d\times \R^d)\}: \|v\|_X<\infty\}$,
where 
\[\|v\|_X = \sup\limits_{0\leq t\leq T}e^{-\lambda t}\|v(t)\|_2 +\Big( \int_0^T\intr e^{-\lambda s}|\Dv v|^2dxd\vi dt\Big)^{\frac{1}{2}},\]
for some $\lambda >0$ to be determined later. We define the map $T:X\rightarrow X$ by $T(w) = u$, where $u$ is the solution to 
\begin{equation}
    \begin{cases}
        \pt u-\lapv u +\vDx u = f- H(\Dv w) \text{ in }(0,T)\times \R^d\times \R^d,\\
        u(0,x,\vi) = g(x,\vi),
    \end{cases}
\end{equation}
where in the above we took the equation forward in time only for notational simplicity. The goal now is to show that $T$ is a contraction on $X$ if $\lambda$ is large enough. But indeed if $C_H>0$ is the Lipschitz constant of $H$, then for $T(w^1)= u^1,T(w^2)= u^2$ by testing against $u^2-u^1$ in the equation of their difference (see the end of this proof on how we justify this), we have 
\[\pt \intr (u^2-u^1)^2(t,x,\vi)dxd\vi +2\intr |\Dv (u^2-u^1)|^2dxd\vi\]
\[\leq C_H\intr |\Dv(w^2-w^1)||u^2-u^1|dxd\vi\leq C_H\epsilon \intr |\Dv(w^2-w^1)|^2 dxd\vi +\frac{C_H}{4\epsilon}\intr |u^2-u^1|^2 dxd\vi    .\]
The above imply 
\[\pt  \Big(e^{-\frac{C_H}{4\epsilon}t}\intr |u^2-u^1|^2dxd\vi   \Big)+2e^{-\frac{C_H}{4\epsilon}t}\intr |\Dv(u^2-u^1)|^2dxd\vi \leq C_H\epsilon e^{-\frac{C_H}{4\epsilon}t} \intr |\Dv(w^2-w^1)|^2dxd\vi, \]
and thus by Gr\"onwall, if we let $\lambda =\frac{C_H}{4\epsilon}$ we have 
\[\|u^2-u^1\|_X^2\leq 4C_H\epsilon \|w^2-w^1\|_X^2.\]
Therefore for $\epsilon>0$ small enough the above map is a contraction in $X$ and thus has a unique fixed point.
Regarding the estimates, we need to test against $u$ in \ref{RenormHJB}. First need to establish integrability for $u$. To this end, note that since $H,f,g\geq 0$, if $w$ is the solution of 
\begin{equation*}
    \begin{cases}
    -\pt w-\lapv w+\vDx w=f \text{ in }(0,T)\times \R^d\times \R^d,\\
    
    w(0,x,\vi)=g(x,\vi) \text{ in }\R^d\times \R^d,
    \end{cases}
\end{equation*}
then by standard comparison we have that 
\[0\leq u(t,x,\vi)\leq w(t,x,\vi) \text{ for all }(t,x,\vi)\in [0,T]\times \R^d\times \R^d.\]
Finally, note that from our assumptions on $f,g$
\[w\in L^p, \text{ for all }p\in [2,\infty]\text{ therefore } u\in L^p \text{ for all }p\in [2,\infty].\]
Now that we may test against $u$ in the equation, the fact that $u\in C([0,T];L^2(\R^d\times \R^d))$ is easy to see due to our assumptions on $f$. The first estimate is obtained by simply testing against $u$ and using the fact that $H$ is Lipschitz with $H(0)=0$. To justify this however, we need to address the integration by parts that occurs. To this end let $\phi:[0,\infty)\rightarrow [0,1]$, be a smooth function such that $\phi(s) = 1$ for $0\leq s\leq 1$ and $\phi(s) =0$ for $s\geq 2$. For $R>0$ we consider the function $\psi_R(x,\vi)= \phi(\frac{|x|^2 + |\vi|^2}{R})$. Testing against $u\psi_R^2$ in equation \ref{RenormHJB} yields
\[-\pt \intr \frac{1}{2}|u|^2 \psi_R^2 dxd\vi +\intr |\Dv u|^2 \psi_R^2 dxd\vi\]
\[+\intr 2u\Dv u \Dv \psi_R \psi_R +2\psi_R \vi\cdot D_x \psi_R u^2 +H(\Dv u)u\psi_R^2 dxd\vi =\intr fu\psi_R^2 dxd\vi. \]
In what follows the constant $C>0$ may change from line to line, however it is independent of $R>0$. We note that
\[\intr \Big|2\psi_R\vi \cdot D_x \psi_R u^2 \Big|=\intr \Big|\frac{2x\cdot \vi}{R}\psi_R \phi'(\frac{|x|^2+|\vi|^2}{R})u^2\Big| \leq \intr \frac{|x|^2+|\vi|^2}{R}\psi_R \phi'(\frac{|x|^2+|\vi|^2}{R})u^2\]
\[\leq \intr 2\psi_R \phi'(\frac{|x|^2+|\vi|^2}{R})u^2\leq C\int_{|x^2|+|\vi|^2\geq R}u^2.\]
Moreover,
\[\intr \Big|2u\Dv u \Dv \psi_R \psi_R \Big|dxd\vi\leq \frac{1}{4}\intr |\Dv u|^2 \psi_R^2dxd\vi +C\int_{|x|^2+|\vi|^2\geq R} u^2 dxd\vi  \]
and 
\[\intr H(\Dv u)u\psi_R^2 dxd\vi \leq \intr \frac{1}{4}|\Dv u|^2 \psi_R^2 +Cu^2\psi_R^2 dxd\vi,\]
\[\intr fu\psi_R^2dxd\vi\leq \|f\|_{2}\|u\|_2. \]
Collecting all the above estimates we have that 
\[-\pt \intr \frac{1}{2}|u|^2 \psi_Rdxd\vi +\intr |\Dv u|^2 \psi_R^2 dxd\vi \leq C(\|f\|_2\|u\|_2+\intr u^2\psi_R^2 dxd\vi +\int_{|x^2|+|\vi|^2\geq R}u^2),\]
recalling that $0\leq u\leq w\in L^2$, the result follows by Gr\"onwall and letting $R\rightarrow \infty$.\\
The second estimates are due to Theorems \ref{Thm1.2FiniteTime} and \ref{Thm2.1FiniteTime} in the Appendix. Finally the $L^{\infty}-$bounds follow by similar arguments as in \cite{degond1986global}, Proposition A.3.\\
\end{proof}

\subsection{Degenerate Fokker-Planck equation}
 All the equations should be understood in the distributional sense, unless stated otherwise. In this subsection we study the following equation 
 \begin{equation}
 \label{eq:ModelEquationForLipschitzFP}
 \begin{cases}
 \pt m-\lapv m-\vDx m-\divv(mb)=0 \text{ in }(0,T)\times \R^d\times \R^d,\\
 m(0,\cdot,\cdot)=m_0(\cdot,\cdot) \text{ in }\R^d\times \R^d.
 \end{cases}
 \end{equation}
 The purpose of this subsection is to show the following theorem:
\begin{theorem}
\label{SolutionToFPequation}
Let $b\in L^{\infty}([0,T]\times \R^d\times \R^d)$ and $m_0$ a density which satisfies \ref{AssumptionInitial}. Then, there exists a unique distributional solution $m\in C([0,T]; L^2(\R^d\times \R^d))$ of equation \eqref{eq:ModelEquationForLipschitzFP}. Furthermore, there exists a $C=C(T,\|b\|_{\infty})>0$, such that 
\[\sup\limits_{t\in [0,T]}\|m(t)\|_2+\|\Dv m\|_{\Lt}+\|D_x^{1/3}m\|_{\Lt}+\|D_t^{1/3}m\|_{\Lt}\leq C\|(1+|\vi|^2)m_0\|_2\]
and a $C_0=C_0(\|b\|_{\infty},T,\|m_0\|_2,\|m_0\|_{\infty})>0$, so that
\[\|m\|_{\infty}\leq C_0.\]
Moreover, $m(t)$ is a probability density for all $t\in (0,T]$. Finally, if $(T-t)\divv(b)\in \Lt$, it follows that
\[[m_t-\vDx m],(T-t)\lapv m\in \Lt.\]
\end{theorem}
The main two assertions in the theorem above are, firstly, the fractional gradient estimates and, secondly, the $L^{\infty}-$bounds. The gradient estimates are the result of Theorem \ref{Thm1.3FiniteTime}, in the appendix. The $L^{\infty}-$bounds can be obtained with a De Giorgi type argument similar to the one found for example in F. Golse, C. Imbert, C. Mouhot and A. Vasseur in \cite{golse2016harnack}, thus we only provide the main steps in Proposition \ref{prop: UpperBoundForm} at the Appendix. For a survey on the De Giorgi type arguments we refer to Mouhot \cite{mouhot2018giorgi}. First a proposition.
\begin{proposition}\label{Prop33}
Assume that $m\in \Lt$, $b\in L^{\infty}\cap \Lt$ and $m_0$, which satisfies \ref{AssumptionInitial}, satisfy equation \eqref{eq:ModelEquationForLipschitzFP} in the distributional sense. Then, $|\vi|^2m,|\vi|^2\Dv m\in \Lt$.
\end{proposition}
\begin{proof}
We may assume that the data are smooth and bounded and obtain the general case by approximation. We test the equation with $|\vi|^4m$ (see Lemma \eqref{Lemma|x|^2+|v|^2m}, on how we may justify this) to obtain
\[\frac{d}{dt}\intr |\vi|^4|m|^2dxd\vi +\intr |\vi|^4|\Dv m|^2dxd\vi \]
\[=-4\intr m|\vi|^2\vi \cdot \Dv mdxd\vi -4\intr |m|^2|\vi|^2 \vi \cdot b dxd\vi -\intr m|\vi|^4 \Dv m\cdot bdxd\vi\]
\[\leq \frac{1}{4}\intr |\vi|^4 |\Dv m|^2dxd\vi +C\intr |m|^2(1+|\vi|^4)dxd\vi \]
\[+4\|b\|_{\infty}\intr |m|^2(1+|\vi|^4)dxd\vi +\intr |m|^2 |\vi|^4dxd\vi +\frac{1}{4}\intr |\Dv m|^2|\vi|^4dxd\vi.\]
It is easy to see that $\sup\limits_{t\in [0,T]}\|m(t)\|_2\leq C\|m_0\|_2,$ therefore the result follows by Gr\"onwall since,
\[\frac{d}{dt}\intr |m|^2 |\vi|^4dxd\vi \leq C\intr |m|^2|\vi|^4dxd\vi +C\|m_0\|_2^2.\]
\end{proof}
\begin{proof} (Theorem \ref{SolutionToFPequation})
Proposition \ref{Prop33}, together with Theorem \ref{Thm1.3FiniteTime}, gives us the result.
\end{proof}

\subsection{Existence of Solutions via the fixed point argument}
In this section we show the main theorem.
\begin{theorem}\label{ExistenceTheoremForWeak}
Let $G,H,F$ and $m_0$ satisfy \ref{AssumptionTerminal},\ref{AssumptionHamilt},\ref{AssumptionCoupling} and \ref{AssumptionInitial}. Then, there exists a unique solution to system \eqref{MFGequation}, according to definition \eqref{WeakSolutionsDef}.
\end{theorem}
\begin{proof}
As mentioned in the introduction, we apply Schauder in the following setting. Let $C_0>0$ be the constant from Theorem \ref{SolutionToFPequation} and consider the closed convex subset $X:=C([0,T];L^2(\R^d\times\R^d ))\cap \{m:\|m\|_{\infty}\leq L\}$ of $C([0,T];L^2(\R^d\times \R^d)_+)$, where $L> 0$, such that $\frac{1}{C_0}\max\{\|m_0\|_{\infty},\|m_0\|_2\}\leq L$. For $\mu\in X$, let $u_{\mu}$ be the solution of
\begin{equation*}
\begin{cases}
-\pt u_{\mu}-\lapv u_{\mu}-\vDx u_{\mu}+H(\Dv u_{\mu})=F(t,x,\vi,\mu(t,x,
\vi)) \text{ in }(0,T)\times \R^d\times \R^d,\\
u_{\mu}(T,x,\vi)=G(\mu(T,x,\vi)) \text{ in }\R^d\times \R^d,
\end{cases}
\end{equation*}
provided by Theorem \ref{EstimatesonHJB}. For this $u_{\mu}$, we then solve
\begin{equation*}
\pt m-\lapv m-\vDx m-\divv(mH_p(\Dv u))=0,\,\,
m(0)=m_0.
\end{equation*}
We set $\Phi(\mu)=m$ which due to the choice of $L$ and the bounds on $m$ implies that $m\in X$. It remains to show that the map is continuous and compact in order to apply Schauder's Fixed Point Theorem.
Continuity is straightforward to check with our given assumptions and will be omitted. For compactness, we proceed as follows. Due to the domain being unbounded we first show that $\lim\limits_{N\rightarrow \infty}\sup\limits_{\mu\in X}\|\Phi(\mu)\mathbf{1}_{B(0,N)^c}\|_2=0,$ where $B(0,N):=\{(t,x,\vi)\in [0,T]\times \R^d\times \R^d:|(x,\vi)|\leq N\}$. This follows directly by the same argument as in Lemma \ref{Lemma|x|^2+|v|^2m}. Furthermore, from Theorem \ref{SolutionToFPequation}, we have
\[\|m\|_2+\|D_{t,x,\vi}^{s}m\|_2\leq C\|m_0\|_2\text{ for some }s>0.\]
Thus, by Kolmogorov–M. Riesz–Fréchet (see for example Brezis \cite{brezis}, Theorem 4.26 and Corollary 4.27) we have compactness of the map. Uniqueness follows from the by-now classical Lasry-Lions monotonicity argument, which we omit.
\end{proof}
We conclude this section with some crucial estimates, which follow directly from the by-now classical Lasry-Lions argument under assumptions \ref{F_LAssum} and \ref{G_LAssum}, so we omit the proof. The computations can be found for example in \cite{porretta2015weak}.
\begin{proposition}\label{EnergyEstimateProp}
Assume that $H:\R^d\rightarrow \R,\,\, F:[0,T]\times \R^d\times \R^d\times \R\rightarrow \R,\,\, m_0:\R^d\times \R^d\rightarrow \R$ and $G:\R^d\times \R^d\times \R\rightarrow \R$ satisfy \ref{AssumptionHamilt},\ref{F_LAssum},\ref{AssumptionInitial} and \ref{G_LAssum}. Let  $(u,m)$ be the weak solution of the MFG system provided by Theorem \ref{ExistenceTheorem}. Then, there exists a constant $C=C(\|m_0\|_1,\|m_0\|_{\infty},T)$, such that
\begin{equation}
\label{EnergyEstimate}
\begin{split}
    \intr G(x,\vi,m(T))dxd\vi +\int_0^T\intr F(x,v,t, m)m dxd\vi ds\\
    +\int_0^T\intr m\Big[ H_p(\Dv u)\cdot \Dv u-H(\Dv u) \Big]dxd\vi \leq C.
\end{split}
\end{equation}
Furthermore, we have the following $L^1$ estimates
\[\sup\limits_{t\in [0,T]}\|u(t)\|_1+\|F(\cdot,m)\|_1+\|F(\cdot,m)m\|_1+\|G(\cdot,m(T))\|_1+\|G(\cdot,m(T))m(T)\|_1\]
\[+\|H(\Dv u)\|_1+\|m |H_p(\Dv u)|^2\|_1\leq C.\]
\end{proposition}
\subsection{Further Regularity of Solutions to the Mean Field Games System, for Lipschitz Hamiltonian}
In this section we study the gain of regularity for solutions to the MFG system \eqref{MFGequation}. In particular, we derive appropriate energy estimates by taking advantage of the coupling.
\begin{theorem}\label{ExtraRegularityLipschitz}
Let $F,G$ satisfy \ref{F_LAssum},\ref{G_LAssum} with constant $c_0$, $H,m_0$ satisfy \ref{AssumptionHamilt},\ref{AssumptionInitial} and $(u,m)$ be a weak solution to system \eqref{WeakSolutionsDef}, according to Definition \eqref{WeakSolutionsDef}. Then, there exists a constant $C=C(c_0,F,G,H,m_0)>0,$ such that
\[\sup\limits_{t\in [0,T]}\|m(t)\|_2+\sup\limits_{t\in [0,T]}\|Dm(t)\|_2+\|D_{\vi,\vi}^2 m\|_2+\|\Dv D_x m\|_2\leq C\]
and 
\[\sup\limits_{t\in [0,T]}\|u(t)\|_2+\sup\limits_{t\in [0,T]}\|Du(t)\|_2+\|D_{\vi,\vi}^2 u\|_2+\|\Dv D_x u\|_2\leq C.\]

\end{theorem}
\begin{proof}
For $i\in \{1,\cdots,d\}$ and $h\in \R\setminus \{0\}$, we denote
\[\delta^h(u)(t,x,\vi):=\frac{u(t,x+he_i,\vi)-u(t,x,\vi)}{h},\delta^h(m)(t,x,\vi):=\frac{m(t,x+he_i,\vi)-m(t,x,\vi)}{h} \]
\[m^h:=m(t,x+he_i,\vi), m^0:=m(t,x+he_i,\vi),\Dv u^h:=\Dv u(t,x+he_i,\vi),\Dv u^0:=\Dv u(t,x,\vi)\]
\[H^h:=H(\Dv u(t,x+he_i,\vi)), H^0:=H(\Dv u(t,x,\vi)),\]
\[F^h:=F(t,x,\vi,m(t,x+he_i,\vi)), F^0:=F(t,x,\vi,m(t,x,\vi)),\]
\[\delta_{x,h}F:=\frac{F(t,x+he_i,\vi,m(t,x+he_i,\vi))-F(t,x,\vi,m(t,x+he_i,\vi))}{h},\]
\[\delta_{x,h}G:=\frac{G(x+he_i,\vi,m(T,x+he_i,\vi))-G(x,\vi,m(T,x+he_i,\vi))}{h}.\]
The equations for $\delta^h u,\delta^h m$ read as follows,
\begin{equation}\label{Equation_uh25}
\begin{cases}
-\pt \delta^h u-\lapv \delta^h u+\vDx \delta^h u+\frac{H^h-H^0}{h}=\frac{F^h-F^0}{h}+\delta_{x,h}F,\\
\delta^h u(T)=\frac{G^h-G^0}{h}+\delta_{x,h}G.
\end{cases}
\end{equation}
\begin{equation}\label{Equation_mh26}
\begin{cases}
\pt \delta^h m-\lapv \delta^h m-\vDx \delta^h m-\divv(\frac{m^h H_p^h-m^0H_p^0}{h})=0,\\
\delta^h m(0)=\delta^h m_0
\end{cases}
\end{equation}
Testing against $\delta^h u$ in \eqref{Equation_mh26}, yields
\[\Bigg[\intr \frac{G^h-G^0}{h}\delta^h m(T)dxd\vi\Bigg]_1+\Bigg[\int_0^T \intr \delta^h m[\frac{F^h-F^0}{h}]dxd\vi dt\Bigg]_2\]
\[+\Bigg[\int_0^T \intr -\delta^h m\frac{H^h-H^0}{h}dxd\vi dt+\int_0^T \intr  \Dv \delta^h u\frac{m^hH_p^h-m^0H_p^0}{h}\Bigg]_3\]
\[=-\int_0^T\intr \delta_{x,h}F\delta^h mdxd\vi dt+\intr \delta^h m_0 \delta^h u(0)dxd\vi -\intr \delta_{x,h}G\delta m^h(T) dxd\vi \]
In the following, we refer to the terms based on the enumeration of the brackets. For the first bracketed term using the monotonicity of $G$ we have
\[\intr \frac{G^h-G^0}{h}\delta^h m(T)dxd\vi=\]
\[\intr \int_0^1 G'(m^0(T)+s(m^h-m^0)(T))ds |\delta^h m|^2(T)dxd\vi \geq c_0\intr |\delta^h m|^2(T)dxd\vi,\]
while for the second term again using the monotonicity of $F$
\[\int_0^T \intr \delta^h m\frac{F^h-F^0}{h}dxd\vi dt=\int_0^T \intr |\delta^h m|^2(t)\int_0^1 F'(m^0(t)+s(m^h-m^0)(t))dsdxd\vi dt\]
\[\geq c_0\int_0^T\intr |\delta^h m|^2(t)dxd\vi dt.\]
We may rewrite the third term as in the proof of uniqueness to see that it is non-negative by the convexity of $H$, indeed it can be written as
\[\int_0^T\intr \frac{m^h}{h^2}\Big[H(\Dv u)-H(\Dv u^h)-H_p(\Dv u^h)\Dv (u-u^h) \Big]\]
\[+\frac{m}{h^2}\Big[H(\Dv u^h)-H(\Dv u)-H_p(\Dv u)\Dv (u^h-u) \Big]dxd\vi dt\geq 0.\]
Continuing, for the right hand side we estimate as follows
\[\delta_{x,h}F=\int_0^1 \partial_{x_i}F(t,x+she_i,\vi,m(t,x+he_i,\vi))hds,\]
hence,
\[\|\delta_{x,h}F\|_2\leq C\|m\|_2,\]
and similarly for $\delta_{x,h}G$. Thus,
\[-\int_0^T\intr \delta_{x,h}F\delta^h mdxd\vi dt+\intr \delta^h m_0 \delta^h u(0)dxd\vi -\intr \delta_{x,h}G\delta m^h(T) dxd\vi\]
\[\leq \frac{c_0}{2}\int_0^T\intr |\delta m^h|^2dxd\vi dt+\frac{c_0}{2}\|\delta m^h(T)\|_2+C\sup\limits_{t\in [0,T]}\|m(t)\|_2^2+\|\delta^h m_0\|_2\|\delta^h u(0)\|_2.\]
Gathering everything together we obtain 
\begin{equation}
\label{d_hm27}
\|\delta^h m(T)\|_2^2+\|\delta^h m\|_2^2\leq C\|\delta^h m_0\|_2\|\delta^h u(0)\|_2.
\end{equation}
We now turn to \eqref{Equation_uh25}. Test, against $\delta^h u$ to obtain 
\[\sup\limits_{t\in [0,T]}\|\delta^h u(t)\|_2+\|\Dv \delta^h u\|_2\leq C(\|\delta^h m(T)\|_2+\|\delta^h m\|_2)\]
and using this estimate in \eqref{d_hm27} provides
\[\|\delta^h m(T)\|_2+\|\delta^h m\|_2\leq C=C(\inf F',\inf G',T,\text{Lip}_H,\text{Lip}_F,\text{Lip}_G,\|D_x m_0\|_2).\]
Testing against $\delta^h m$ in \eqref{Equation_mh26} yields
\[\sup\limits_{t\in [0,T]}\|\delta^h m(t)\|_2+\|\Dv \delta^h m\|_2\leq C(\|\delta^h m_0\|_2+\|\Dv \delta^h u\|_2)\leq C.\]
Since the bounds are independent of $h$, we have shown that 
\[\sup\limits_{t\in [0,T]}\|D_x m(t)\|_2+\sup\limits_{t\in [0,T]}\|D_x u(t)\|_2+\|\Dv D_x u\|_2+\|\Dv D_x m\|_2\leq C.\]
Now, using these bounds, we repeat the process for the derivatives with respect to $\vi$. We use completely symmetric notation as in the above case, for example $\delta_{\vi}^h u:=\frac{u(t,x,\vi+he_i)-u(t,x,\vi)}{h}$. The equations satisfied by $\delta_{\vi}^h u,\delta_{\vi}^h m$ are similar with the exception of the $\vDx $ term. They read
\begin{equation*}
\begin{cases}
-\pt \delta_{\vi}^h u-\lapv \delta_{\vi}^h u+e_{\vi,i}D_x u^h+\vDx \delta_{\vi}^h u+\frac{H^h-H^0}{h}=\frac{F^h-F^0}{h}+\delta_{\vi,h}F,\\
u^h(T)=\frac{G^h-G^0}{h}+\delta_{\vi,h}G
\end{cases}
\end{equation*}
and
\begin{equation*}
\begin{cases}
\pt \delta_{\vi}^h m-\lapv \delta_{\vi}^h m-e_{\vi,i}D_x m^h-\vDx \delta_{\vi}^h m-\divv(m^h\frac{H_p^h-H_p^0}{h}+\delta^h mH_p^0)=0,\\
\delta_{\vi}^h m^0=\delta_{\vi}^h m_0.
\end{cases}
\end{equation*}
The argument is completely symmetrical with the only difference being the presence of $D_x u^h,D_x m^h$. However, these terms are bounded from the previous case. We thus obtain bounds of the form
\[\sup\limits_{t\in [0,T]}\|\Dv m(t)\|_2+\sup\limits_{t\in [0,T]}\|\Dv u(t)\|_2+\|D_{\vi,\vi}^2 u\|_2+\|D_{\vi,\vi}^2 m\|_2\leq C.\]
\end{proof}

\section{Quadratic Hamiltonian}
In this section we will show existence and uniqueness for renormalized solutions to the MFG system. All the ideas and proofs in this section are entirely motivated or even parallel to the original work in \cite{porretta2015weak}.\par
To motivate some of the technical steps we outline the strategy. The plan is to approximate a given Hamiltonian $H$ with quadratic growth by a sequence of Lipschitz Hamiltonians $H^{\epsilon}$(see bellow for definition), for which we have shown the existence of solutions $(\ue,\me)$ in the previous section and show that these solutions converge to a renormalized solution. A crucial structural estimate, as pointed out in \cite{porretta2015weak}, is that $\sup\limits_{\epsilon}\|\me |H_p^{\epsilon}(\Dv \ue)|^2\|_1<\infty$, which is shown in Proposition \ref{EnergyEstimateProp}. This estimate, along with $L^2-$bounds on $\Dv \ue$, allows us to conclude the convergence (up to a subsequence) to a renormalized solution of $\{\me\}_{\epsilon}$. The bounds for the HJ equation are straightforward and mostly follow the classical techniques of the non-degenerate case, with the exception of the $L^1-$compactness for the $\ue$ which is due to Theorem \ref{L1-comp} in \cite{diperna1988fokker} and the technical Lemma in the Appendix.\par
In the rest of the paper we consider a fixed Hamiltonian $H$ that satisfies \ref{Assumption2onH}. Furthermore, following \cite{porretta2015weak}, we consider the following Lipschitz approximations
\begin{equation}
\label{Hepsilon}
H^{\epsilon}(p):=\frac{H(p)}{1+\epsilon H^{\frac{1}{2}}(p)} \text{ for }\epsilon>0.
\end{equation}
The following are shown in \cite{porretta2015weak}
\begin{proposition}
\label{BoundsHepsilon}
The functions $H^{\epsilon}$ are Lipschitz in $p$ and satisfy 
\[H_p^{\epsilon}\cdot p-H^{\epsilon}(p)\geq cH^{\epsilon}(p), |H_p^{\epsilon}|^2\leq CH^{\epsilon},\]
for some constants $c>0$,$C>0$ independent of $\epsilon$.
\end{proposition}

\subsection{Analysis of Degenerate Fokker-Planck equation}
In this subsection, we study the following Fokker-Planck equation
\begin{equation}
\label{FPEquationQuadratic}
\begin{cases}
\pt m-\lapv m-\vDx m-\divv(mb)=0 \text{ in }(0,T)\times \R^d\times \R^d,\\
m(0,x,\vi)=m_0(x,\vi) \text{ in }\R^d\times \R^d,
\end{cases}
\end{equation}

 Our approach is a parallel of the techniques from \cite{boccardo1997nonlinear} in the Hypoelliptic case.
\begin{definition}\label{WeakFPSolQuadr}
We say that $m$ is a weak solution of \eqref{FPEquationQuadratic}, if $m\in L^1\cap L^{\infty}([0,T]\times \R^d\times \R^d),$ with $\Dv m\in L^2([0,T]\times \R^d\times \R^d)$, $m_0$ satisfies \ref{AssumptionInitial}, $m|b|^2\in L^1([0,T]\times \R^d\times \R^d)$, and \eqref{FPEquationQuadratic} is satisfied in the distributional sense.
\end{definition}
\begin{lemma}\label{Lemma|x|^2+|v|^2m}
Let $(m,b,m_0)$ be a weak solution of \eqref{FPEquationQuadratic} according to definition \ref{WeakFPSolQuadr}. Then, there exists a constant $C=C(d,T,\|m|b|^2\|_1,\|(1+|x|^2+|\vi|^2)m_0\|_1)$, such that for all $t\in [0,T]$
\[\intr (|x|^2+|\vi|^2+1)m(t,x,\vi)dxd\vi \leq C.\]
\end{lemma}
\begin{proof}
Formally the result follows immediately by testing against $(|x|^2+|\vi|^2)$ and applying standard methods. However, this needs to be justified given that $(|x|^2+|\vi|^2)$ is unbounded. This requires some technical steps  which we present in detail, hence the lengthy computations. First assume that $b,m_0$ are smooth and compactly supported. For $R>0$ consider a bump function $\psi_R:\R^d\times \R^d\rightarrow [0,1]$, such that $\psi_R\Big|_{B(0,R)}\equiv 1$ and $\text{spt}(\psi_R)\subset B(0,R+1)$. Fix a $t_0\in [0,T]$ and let $\phi_R:[0,t_0]\times \R^d\times \R^d\rightarrow \R$ be the smooth solution of the adjoint equation (see for example E. Priola \cite{Adjoint}, Theorem 5.3)
\begin{equation}
    \label{AdjointEq}
    \begin{cases}
    -\pt \phi_R-\lapv \phi_R+\vDx \phi_R+b\cdot \Dv \phi_R =0 \text{ on }[0,t_0)\times \R^d\times \R^d,\\
    \phi_R(t_0,x,\vi)=(|x|^2+|\vi|^2)\psi_R(x,\vi) \text{ on }\R^d\times \R^d.
    \end{cases}
\end{equation}
A priori, $\phi_R$ is bounded by a constant depending only on $R,b,T$. We claim that there exists a constant $C>0$ independent of $R>0$, such that 
\[\phi_R(t,x,\vi)\leq C(1+|x|^2+|\vi|^2) \text{ for all }(t,x,\vi)\in [0,t_0]\times \R^d\times \R^d.\]
Indeed, for $A,B>0$ large enough to be determined later, let $w(t,x,\vi)=Ce^{-At}(1+|x|^2+|\vi|^2)-B(t-t_0)$, which satisfies
\[ -\pt w-\lapv w+\vDx w+b\cdot \Dv w=Ce^{-At}(A(1+|x|^2+|\vi|^2)-2d+2\vi \cdot x+b\cdot \vi)+B\]
\[\geq (B-2dCe^{-At}-\|b\|_{\infty}^2)+Ce^{-At}(A-\frac{3}{2})(|x|^2+|\vi|^2)\geq 1,\]
if $A,B>0$ are large enough. In particular let $A=2$ and for any choice of $C>0$ we set $B=1+2dCe^{-2t}-\|b\|_{\infty}^2$, so that the above inequality is satisfied. Furthermore, at $t=t_0$ we have that
\[w(t_0,x,\vi)=Ce^{-2t_{0}}(|x|^2+|\vi|^2)\geq (|x|^2+|\vi|^2)\psi_R(x,\vi)= \phi_R(t_0,x,\vi) \text{ for all }(x,\vi)\in \R^d\times \R^d,\]
if say $C>e^{2t_0}$, in particular however $C$ can be chosen independent of $R>0$. Finally, for each $R>0$ fixed, the function 
\[E(t,x,\vi)=w-\phi_R\]
is coercive in $(x,\vi)$, that is for each fixed $t\in [0,t_0]$,
\[\lim\limits_{|(x,\vi)|\rightarrow \infty}E(t,x,\vi)=\infty.\]
Thus by classical arguments we find that the minimum of $E$ is achieved at $t=t_0$, which shows the claim. To conclude the proof of the Lemma, we test against $\phi_R$ in equation \eqref{FPEquationQuadratic}, which yields
\[\intr m(t_0) (|x|^2+|\vi|^2)\psi_R(x,\vi)dxd\vi =\intr \phi_R(0,x,\vi)m_0(x,\vi)dxd\vi\]
\[\leq C\intr m_0(|x|^2+|\vi|^2+1)dxd\vi=C\|m_0(1+|x|^2+|\vi|^2)\|_1.\]
With the above bounds we may now test equation \eqref{FPEquationQuadratic} against $(1+|x|^2+|y|^2)$, which yields
\[\pt \intr (1+|x|^2+|\vi|^2)m(t)dxd\vi =\intr 2dm(t) -2 x\cdot \vi m(t) +2m(t)\vi\cdot  b dxd\vi\]
\[\leq 2d m(t) +\int (1+|x|^2+|\vi|^2)m(t)+m|\vi|^2+m|b|^2dxd\vi\leq (2d+2)\intr (1+|x|^2+|\vi|^2)m(t)dxd\vi +\|m|b|^2\|_1, \]
and so by Gr\"onwall we obtain that for some constant $C=C(d,T,\|m|b|^2\|_1,\|(1+|x|^2+|\vi|^2)m_0\|_1)>0$ 
\[\|(|x|^2+|\vi|^2+1)m_0\|_1\leq C.\]
The general case follows by approximation with smooth data.
\end{proof}
In the following Proposition we will need the following estimate, which may be found in in Folland \cite{folland1975subelliptic}.
\begin{proposition}\label{Folland}[\cite{folland1975subelliptic}, Theorem 5.14] Let $\Gamma$ denote the fundamental solution of the operator $\pt -\lapv-\vDx$ in the space $\R^d\times \R^d$. Assume that $p,q\in (1,\infty)$ are such that 
\[\frac{1}{p}=\frac{1}{q}-\frac{1}{Q+2},\]
where $Q=d+2$. For a function $f\in L^q$ we define
\[g(t,x,\vi):=\int_0^T\intr \Dv \Gamma(t-s,x,\vi,y,w)f(s,w,y)dydwds.\]
Then, there exists a constant $C=C(p,q,d)$ such that 
\[\|g\|_p\leq C\|f\|_q.\]
\end{proposition}
 \begin{proposition}\label{GainIntegrability_m}
Let $(m,b,m_0)$ be a weak solution of \eqref{FPEquationQuadratic}, according to definition \ref{WeakFPSolQuadr}. Then, there exists a dimensional constant $C=C(d)>0$ and a constant $C_0=C_0(m_0)>0$, such that 
\[\|m|b|^2\|_{\frac{d+4}{d+3}}+\|m\|_{\frac{d+4}{d+2}}\leq C\|m|b|^2\|_1+C_0.\]
\end{proposition}
\begin{proof}
 Let $\Gamma$ denote the fundamental solution of the operator $\pt-\lapv -\vDx $. From the equation satisfied by $m$ we obtain
\[m(x,\vi,t)=-\int_0^t\intr D_{\vi}\Gamma (t-s,x,\vi,y,w) mb(s,w,y)dydw ds+C(m_0)(t,x,\vi)\]
where
\[C(m_0)(t,x,\vi)=\intr \Gamma(t,x,\vi,y,w)m_0(y,w)dy dw.\]
From Proposition \eqref{Folland} above, we have that
\[\|m\|_p\leq C\|mb\|_q\]
where
\[\frac{1}{p}=\frac{1}{q}-\frac{1}{Q+2},\]
and $Q=d+2$. Moreover, by H\"older
\[\int_0^T\intr |m|^q|b|^qdxd\vi dt\]
\[\leq \Big(\int_0^T\intr |m|^{\frac{q}{2-q}}dxd\vi dt\Big)^{\frac{2-q}{2}}\Big(\int_0^T\intr m|b|^2dxd\vi \Big)^{\frac{q}{2}}=C\|m\|_{\frac{q}{2-q}}^{\frac{q}{2}}.\]
Hence, we can have a gain of integrability if we require that
\[p=\frac{q}{2-q}\iff \frac{2-q}{q}=\frac{1}{q}-\frac{1}{Q+2}\iff \frac{1}{q}-1=-\frac{1}{Q+2}\iff \frac{1}{q}=\frac{Q+1}{Q+2},\]
therefore
\[q=\frac{Q+2}{Q+1} \text{ and }p=\frac{Q+2}{2Q+2-Q-2}=\frac{Q+2}{Q}.\]
\end{proof}
\begin{proposition}\label{Log-EstimatesPropForm}
Let $(m,b,m_0)$ be a weak solution of \eqref{FPEquationQuadratic} according to Definition \ref{WeakFPSolQuadr}, with $b\in L^2([0,T]\times \R^d\times \R^d;\R^d)$. Then, there exists a constant $C(\|m_0\log(m_0)\|_1,\|m|b|^2\|_1)>0$ such that
\[\sup\limits_{t\in [0,T]}\|m(t)\log(m(t))\|_1+\|\Dv (\sqrt{m})\|_2\leq C.\]

\end{proposition}
\begin{proof}
 For $\delta>0$, define $w(x)=\log(x+\delta)$ and $W(x)=(x+\delta)\log(x+\delta)-\delta\log(\delta)$. Test against $w(m)$ in \eqref{FPEquationQuadratic} ($m\in L^{\infty}\cap L^1$ and so $w(m)\in L^{\infty}, W(m)\in L^1$) to obtain that for each $t\in [0,T]$
\[\intr W(m(t))dxd\vi +\int_0^t\intr \frac{|\Dv m|^2}{(m+\delta)}dxd\vi ds =-\int_0^t\intr \frac{m}{m+\delta}\Dv m\cdot bdxd\vi dt\]
\[+\intr W(m_0)dxd\vi \]
\[\leq \frac{1}{2}\int_0^T\intr \frac{|\Dv m|^2}{(m+\delta)}dxd\vi dt+\frac{1}{2}\|m|b|^2\|_1+\intr W(m_0)dxd\vi .\]
Letting $\delta\rightarrow 0$ yields
\[\int m(t)\log(m(t))dxd\vi +\int_0^T\intr \frac{|\Dv m|^2}{m}dx\vi ds\leq C(\|m|b|^2\|_1+\|m_0\log(m_0)\|_1)\]
where $C>0$ is a universal constant. It remains to show that $m(t)\log(m(t)\in L^1$. This is shown for example in \cite{diperna1988fokker}, under the conditions
\begin{enumerate}
    \item\label{NormXbound}$\|m(t)(1+|x|^2+|\vi|^2)\|_1<\infty$
    \item\label{mlogm} $\intr m(t)\log(m(t))<\infty.$
\end{enumerate}
Condition \ref{NormXbound} follows from Lemma \ref{Lemma|x|^2+|v|^2m}, while condition \ref{mlogm} is shown above.
\end{proof}

We now proceed with gradient estimates for the measure.
\begin{theorem}\label{GradientEstimatesForm} Let $(m,b,m_0)$ be a weak solution of \eqref{FPEquationQuadratic} according to Definition \ref{WeakFPSolQuadr}. Then, there exist $s\in (0,1)$ and $q\in (1,\infty)$, such that 
\[\|D^sm\|_q\leq C,\]
where $C$ depends only on $m_0,d,T,\|m|b|^2\|_1$ and in particular not on $\|\Dv m\|_2$.
\end{theorem}
\begin{proof}
The constant $C>0$ that appears in this proof is subject to change from line to line and depends only on $m_0,d,T$. The technique that follows is the same as in \cite{boccardo1997nonlinear}.
 In the original equation \eqref{FPEquationQuadratic} we test against $\phi(m)$ for $\phi(s)=s$ for $s\in [0,1]$ and $\phi(s)=1, s\geq 1, \phi(s)=0, s\leq 0$. This yields
\[\intr \Phi(m(t))dxd\vi +\int_0^T\intr \phi'(m)|\Dv m|^2=\]
\[-\int_0^T\intr \phi'(m)\Dv m H_p mdxd\vi dt+\intr \Phi(m_0)dxd\vi \]
\[\leq \frac{1}{2}\int_0^T\intr \phi'(m)|\Dv m|^2dxd\vi +\int_{|m|\leq 1} |m|^2|b|^2dxd\vi+C(m_0) .\]
Since $|m|^2\leq |m|$ on $|m|\leq 1$, we obtain 
\[\int_{\{|m|\leq 1\}}|\Dv m|^2dxd\vi \leq C.\]
For $k\in \N$ we define $\phi_k$ by
\begin{equation}
\phi_k(s):=
\begin{cases}
0, s\leq k-1,\\
s-(k-1), s\leq k,\\
1, s\geq 1,
\end{cases}
\end{equation}
and $\Phi_k(t):=\int_0^T \phi_k(s)ds$. Testing against $\phi_k(m)$ in the equation yields
\begin{equation}
\label{Phi_kTestA_k}
\begin{split}
    \int_0^T\intr \Phi_k(m(T))+\int_0^T\intr \phi_k'(m)|\Dv m|^2dxd\vi dt\\
    =-\int_0^T\intr\phi_k' \Dv m b mdxd\vi dt+\intr \Phi_k(m_0)dxd\vi .
\end{split}
\end{equation}
Additionally
\[\intr \Phi_k(m_0)dxd\vi\leq \|m_0\|_2+\|m_0\|_1\leq C\]
and 
\[0\leq \int_0^T\intr \Phi_k(m(T)).\]
For $A_k:=\{k-1\leq |m|\leq k\}, k\in \N$, equation \eqref{Phi_kTestA_k} yields
\[\int_{A_k}|\Dv m|^2dxd\vi dt\leq \frac{1}{2k}\int_{A_k}m|\Dv m|^2dxd\vi + Ck\int_{A_k}m|b|^2dxd\vi+C, \text{ for all }k\in \N. \]
Moreover,
\[\int_{A_k}m|\Dv m|^2dxd\vi dt\leq k\int_{A_k}|\Dv m|^2dxd\vi dt\]
hence, by summing for $k=2,\cdots,$ for $\lambda>1$, we obtain
\[\int_{|m|\geq 1 }\frac{|\Dv m|^2}{(1+m)^{\lambda}}dxd\vi dt\leq \sum\limits_{k=1}^{\infty} \frac{k}{(1+k)^{\lambda}}\int_{A_k}m|b|^2dxd\vi dt+\frac{C}{k^{\lambda}}<\infty.\]
Thus,
\[\int_{m>1} |\Dv m|^qdxd\vi \leq\Big[\int_{m>1}\frac{|\Dv m|^2}{(1+m)^{\lambda}} \Big]^{q/2}\Big[\int_{m>1} (1+m)^{\frac{\lambda q}{2-q}}dxd\vi \Big]^{\frac{2-q}{2}}.\]
Next, using that 
\[(a+b)^{\lambda}\leq 2^{\lambda}\max\{a^{\lambda},b^{\lambda}\}\leq C(a^{\lambda}+b^{\lambda})\]
and 
\[|\{|m|>1\}|\leq \|m\|_1=1,\]
we obtain
\[\int_{|m|>1}(1+m)^{\frac{\lambda q}{2-q}}dxd\vi\leq C\Big(|\{m>1\}|^{\frac{\lambda q}{2-q}}+\int_{\R^d\times \R^d}|m|^{\frac{\lambda q}{2-q}}dxd\vi  \Big)\leq C(1+\int_{\R^d\times \R^d}|m|^{\frac{\lambda q}{2-q}}dxd\vi).\]
Hence,
\begin{equation}\label{FirstStepGradient-m}
    \int_{m>1} |\Dv m|^qdxd\vi\leq \Big[\int_{m>1}\frac{|\Dv m|^2}{(1+m)^{\lambda}} \Big]^{q/2} \Big(1+\int_{\R^d\times \R^d}|m|^{\frac{\lambda q}{2-q}}dxd\vi\Big)^{\frac{2-q}{2}}.
\end{equation}
Integrate in time inequality \eqref{FirstStepGradient-m}, and apply H\"olders inequality for $\frac{2}{q},\frac{2}{2-q},$ to obtain for some $C=C(T,\lambda,q,\|\frac{\Dv m}{(1+m)^{\frac{\lambda}{2}}}\mathbf{1}_{m\leq 1}\|_2)>0$
\[\int_{m>1} \|\Dv m(t)\|_{q}^qdxd\vi dt\leq \Big( \int_{m>1}\frac{|\Dv m|^2}{(1+m)^{\lambda}}dxd\vi dt \Big)^{\frac{q}{2}}\Big(1+\int_0^T\|m(t)\|_{\frac{\lambda q}{2-q}dt}^{\frac{\lambda q}{2-q}} dt\Big)^{\frac{2-q}{q}}\]
\[\leq C(1+\Big(\int_0^T\|m(t)\|_{\frac{\lambda q}{2-q}}^{\frac{\lambda q}{2-q}}dt\Big)^{\frac{2-q}{2}} )\]
The Fractional Gagliardo-Niremberg inequality gives us
\[\|m(t)\|_{\sigma}\leq C\|D^s m\|_{q}^{\theta}\|m(t)\|_{\rho}^{1-\theta},\]
where
\begin{equation}
    \label{GagliardoNirembergExponents}
    \frac{1}{\sigma}=\theta(\frac{1}{q}-\frac{s}{n})+\frac{1-\theta}{\rho},
\end{equation}
and $C=C(s,q,n,\theta,\rho)>0$, we refer for example to \cite{brezis2018gagliardo}. We can easily obtain the following time dependent version,
\[\int_0^T \|m(t)\|_{\sigma}^{\sigma}dt\leq C\sup\limits_{t}\|m(t)\|_{1}^{\sigma(1-\theta)}\int_0^T \|D^s m\|_q^{\theta \sigma}dt\leq C\int_0^T \|D^s m\|_q^{\theta \sigma}dt.\]
Set
\[\theta=\frac{q}{\sigma},\rho=1, \sigma=\frac{\lambda q}{2-q},\]
which implies that
\[\frac{1}{\sigma}=\frac{q}{\sigma}(\frac{1}{q}-\frac{s}{n})+1-\frac{q}{\sigma}=\frac{1}{\sigma}-\frac{q s}{\sigma n}+1-\frac{q}{\sigma }\]
thus,
\[\frac{q s}{\sigma n}=1-\frac{q}{\sigma }\implies \sigma=\frac{q s}{n}+q\implies\sigma =q(\frac{ s}{n}+1)\]
and so
\[q(\frac{s}{n}+1)=\frac{\lambda q}{2-q}\implies \lambda =(2-q)(1+\frac{s}{n})\]
which is a valid choice as long as 
\[(2-q)(1+\frac{s}{n})>1\implies q<2-\frac{n}{n+s}\]
thus our restrictions on $q$ is that 
\[1<q<2-\frac{n}{n+s}.\]
Continuing with the above analysis for the above choices of parameters we obtain 
\[\int_{m>1} \|\Dv m(t)\|_{q}^qdxd\vi dt\leq C(1+\Big(\int_0^T\|m(t)\|_{\sigma}^{\sigma}dt\Big)^{\frac{2-q}{2}} )\leq C\Big(1+\int_0^T \|D^s m\|_q^qdt\Big)^{\frac{2-q}{2}}. \]
Therefore for some $\alpha\in (0,1)$
\[\|\Dv m\|_q\leq C(\|\Dv m\mathbf{1}_{m\leq 1}\|_q+\|\Dv m \mathbf{1}_{m>1}\|_q)\]
\[\leq C(1+\|\Dv m\|_1^{\alpha}\|\Dv m \mathbf{1}_{m\leq 1}\|_2^{1-\alpha}+\|D^s m\|_q^{\frac{2-q}{2}}),\]
and by using the estimate from Proposition \ref{Log-EstimatesPropForm}, we obtain 
\[\|\Dv m\|_1=\|\sqrt{m}\Dv \sqrt{m} \|_1\leq \|\Dv \sqrt{m}\|_2,\]
therefore
\[\|\Dv m\|_q\leq C(1+\|D^s m\|_q^{\frac{2-q}{2}}).\]
By Theorem \ref{Thm2.1FiniteTime}, we have that 

\[\|D_x^s m\|_q\leq C(1+\|\Dv m\|_q+\|m|b|^2\|_q+\|m\|_q)\]
\[\leq C(1+\|D^s m\|_q^{\frac{2-q}{2}})\]
Thus by choosing $q$ so that $\|m|b|^2\|_q+\|m\|_q\leq C$ from Proposition \ref{GainIntegrability_m}, the result follows.

\end{proof}
\begin{theorem}\label{CompactnesL^1Form}
Let $\{(m^n,b^n,m_0)\}_{n\in \N}$ be a sequence of weak solutions to \eqref{FPEquationQuadratic} according to definition \ref{WeakFPSolQuadr}, such that 
\[\sup\limits_{n\in \N}\Big(\|m^n|b^n|^2\|_{1}+\|b^n\|_2\Big)<\infty.\]
Then, the set $\{m^n\}_{n\in \N}$ is compact in $L^1([0,T]\times \R^d\times \R^d)$.
\end{theorem}
\begin{proof}
From Proposition \ref{GradientEstimatesForm}, we have that 
\[\|m^n\|_r+\|D^s m^n\|_q\leq C \text{ for all }n\in\N \text{ and some }r>1, s\in (0,1).\]
The result about the compactness in $L^1([0,T]\times L^1(\R^d\times \R^d))$ now follows by the results in \cite{simon1986compact}, with a slight modification due to the unbounded domain. We sketch the argument. For $R>0$, let $\phi_R(x,\vi):=\psi_R(x)\psi_R(\vi)$, where $\psi_R$ are standard non-negative cutoff functions with support in $B(0,R)\subset \R^d$. The, equation satisfied by $m^R:=m\phi_R$, reads
\[\pt m^R-\lapv m^R-\vDx m^R-\divv(m^R b)=\Dv \phi^R m b-m\lapv (\phi^R)-2\Dv \phi^R\Dv m-m\vDx \phi^R.\]
Next for $\frac{1}{p}+\frac{1}{q}=1$, we set $X:=W^{s,q}(B_{\R^d\times \R^d}(0,R)), B:=L^q(B_{\R^d\times \R^d}(0,R))$ and $Y:=W^{-1,p}(B_{\R^d\times \R^d}(0,R))$. Space $X$ embeds compactly in $B$ and $B$ embeds continuously in $Y$. Since $m_n^R$ are bounded in $L^q(0,T,X)$ and $\pt m_n^R$ is bounded in $L^q(0,T,Y)\subset L^1((0,T),Y)$. Therefore from Corollary 4 in \cite{simon1986compact}, for each fixed $R>0$ the sequence $m_n^R$ is compact in $L^q(0,T,B)=L^q(0,T,B_{\R^d\times \R^d}(0,R))\subset L^1(0,T,B_{\R^d\times \R^d}(0,R))$. Combining the above with the estimate $\sup\limits_{n,t}\int_{B(0,R)^c}m^n(t,x,\vi)dxd\vi \rightarrow 0 \text{ as }R\rightarrow \infty,$
from Lemma \ref{Lemma|x|^2+|v|^2m}, yields the strong convergence in $L^1([0,T]\times \R^d\times \R^d)$.
\end{proof}
\begin{proposition}\label{WeakLimitFormQuadratic}
Let $\{(m^n,b^n,m_0)\}_{n\in \N}$ be a sequence of weak solutions to \eqref{FPEquationQuadratic} according to definition \ref{WeakFPSolQuadr}, such that 
\[\sup\limits_{n\in \N}\Big(\|m^n|b^n|^2\|_{1}+\|b^n\|_2\Big)<\infty \]
and
\[b^n\rightarrow b \text{ almost everywhere, for some }b\in L^2([0,T]\times \R^d\times \R^d).\]
Then, there exists a $m\in L^1([0,T]\times \R^d\times \R^d)$, such that up a subsequence $m^n\rightarrow m, m^nb^n\rightarrow mb$ in $L^1([0,T]\times \R^d\times \R^d)$. Furthermore, the set $\{m^n\}_{n\in \N}$ is compact in $C([0,T]; \mathcal{P}_1(\R^d\times \R^d))$. Finally, $m$ is a distributional solution of \eqref{FPEquationQuadratic}.
\end{proposition}
\begin{proof}
From Theorem \ref{CompactnesL^1Form}, there exists an $m\in L^1([0,T]\times \R^d\times \R^d)$ and a subsequence(still denoted by $\{m_n\}_{n\in \N}$) such that $\|m_n-m\|_1\rightarrow 0$. Furthermore, from Lemma \ref{Lemma|x|^2+|v|^2m} we have that
\[\limsup\limits_{R\rightarrow \infty}\sup\limits_{n\in \N}\int_0^T\int_{B_R^c}|m^n||b^n|dxd\vi\]
\[\leq\limsup\limits_{R\rightarrow \infty}\sup\limits_{n\in \N} \Big( \int_0^T\int_{B_R^c}|m^n|dxd\vi dt \Big)^{\frac{1}{2}}\Big(\int_0^T\int |m^n||b^n|^2dxd\vi dt \Big)^{\frac{1}{2}}=0.\]
The above combined with Proposition \ref{GainIntegrability_m} yields that the sequence $\{m^nb^n\}_{n\in \N}$ is uniformly integrable, which together with the almost everywhere convergence gives us that the limit $m$ is in fact a distributional solution of \eqref{FPEquationQuadratic}.\\ 
Next, we show the claim about the compactness in $C([0,T];\mathcal{P}_1(\R^d\times \R^d))$. From Lemma \ref{Lemma|x|^2+|v|^2m} the set $\{m^n(t)\}_{n\in \N}$ is compact in $\mathcal{P}_1(\R^d\times \R^d)$ for each $t\in [0,T]$. The result about compactness in $C([0,T];\mathcal{P}_1(\R^d\times \R^d))$, will follow once we obtain H\"older time continuity. However this follows by typical arguments such as the one found in the notes of Cardaliaguet \cite{cardaliaguet2010notes}.
\end{proof}

\begin{theorem}\label{RenormalizedSolutionForm}
Let $\{(m^n,b^n,m_0)\}_{n\in \N}$ be a sequence of weak solutions to \eqref{FPEquationQuadratic} according to definition \ref{WeakFPSolQuadr}.
Assume furthermore that $\sup\limits_n \|b^n\|_2<\infty,$
and that the assumptions of Proposition \ref{WeakLimitFormQuadratic} are satisfied. Then, the limit $m$ provided by Proposition \ref{WeakLimitFormQuadratic} is a renormalized solution according to Definition \ref{RenormFPDefinition}.
\end{theorem}
\begin{proof}
Let $S:\R\rightarrow \R$, such that 
$S\in W^{1,\infty}(\R) \text{ and that }S' \text{ has compact support.}$ Then, for each $n\in \N$ we have
\begin{equation}\label{RenormS}
    \pt S(m^n)-\lapv S(m^n)-\vDx S(m^n)-\divv(S'(m^n)m^nb^n)+S''(m^n)\Dv m^n m^n b^n+S''(m^n)|\Dv m^n|^2=0.
\end{equation}
Since $\{m^n|b^n|^2\}_{n\in \N}$ is uniformly bounded in $L^1([0,T]\times \R^d\times \R^d)$, we obtain that 
\[\lim\limits_{k\rightarrow \infty}\sup\limits_{n\in \N}\frac{1}{k}\int_{k<m^n<2k}|\Dv m^n|^2dxd\vi ds=0,\]
just as in Theorem 6.1 of \cite{porretta2015weak}. It remains to show that for a fixed $k\in \N$, we have the following convergence
$\Dv (m^n\wedge k)\rightarrow \Dv (m\wedge k) \text{ strongly in }L^2.$ To show the strong convergence of the truncations, it is enough to show that
\[\|\Dv \log(1+m_n)-\Dv \log(1+m)\|_{\Lt}\rightarrow 0.\]
 The argument that follows is entirely due to DiPerna-Lions in \cite{diperna1988fokker}. We only present some of the main estimates since we have a slightly different setup. We look at $g^n=\log(1+m_n)$ and the corresponding equation they satisfy. From Proposition \ref{Log-EstimatesPropForm} we have that $\sup\limits_{n\in \N}\|\Dv g^n\|_2<\infty$ and so without loss of generality we may assume that $\Dv g^n$ converges weakly in $L^2$ to $\Dv g$, where $g=\log(1+m)$. Therefore, there exists a non-negative bounded measure $\mu $ (in the sense that $\int_0^T\intr d\mu<\infty$) on $(0,T)\times \R^d\times \R^d$ such that 
\[|\Dv g^n|^2\rightarrow |\Dv g|^2+\mu\]
in the distributional sense. It remains to show that $\mu$ is identically zero. First, for each $n\in \N$ we let $\beta=\log(1+ t)$ and $g^n=\beta(m^n).$ The functions $g^n$ satisfy
\[\pt g^n-\lapv g^n-\vDx g^n-\divv(\frac{m^n}{1+m^n}b^n)=|\Dv g^n|^2+\frac{m^n}{1+m^n}b^n\Dv g^n\]
\[g^n(0)=\log(1+m_0).\]
Again, just as in \cite{diperna1988fokker}, we set $\Phi_{s,R}^n(t)=\exp(s t\wedge R)$ and $\Psi_{s,R}^n(t):=\int_0^T \Phi_{s,R}^n(\theta)d\theta,$
for some $0<s<1$. Test the equation against $\Phi_{s,R}^n(g^n) \phi$, where $\phi \in C_c((0,T))$, which yields
\[-\int_0^T\intr \Psi_{s,R}^n(g^n)\phi'(t)dxd\vi dt\]
\[+\int_0^T\intr s\phi |\Dv g^n|^2\mathbf{1}_{g^n\leq R}\Phi_{s,R}^n(g^n)+s\Phi_{s,R}^n(g^n)\mathbf{1}_{g^n\leq R}\Dv g^n \frac{m^n}{1+m^n}b^ndxd\vi dt\]
\[=\int_0^T\intr \Phi_{s,R}^n(g^n)\phi |\Dv g^n|^2+\phi\Phi_{s,R}^n(g^n)\frac{m^n}{1+m^n}b^n\Dv g^n,\]
or equivalently
\begin{equation}
    \label{EqTestAgainPhin}
    \begin{split}
    &-\int_0^T\intr \Psi_{s,R}^n(g^n)\phi'(t)dxd\vi dt=\\
    &\int_0^T\intr \phi \Phi_{s,R}^n(g^n)\Big[\Big(|\Dv g^n|^2-s|\Dv g^n|^2\mathbf{1}_{g^n\leq R} \Big)+\frac{m^n \cdot b^n}{1+m^n}\Big( \Dv g^n-s\Dv g^n\mathbf{1}_{g^n\leq R} \Big) \Big]dxd\vi dt
    \end{split}
\end{equation}
\[=(I)+(II).\]
Now we bound each term,
\[|(I)|\leq \|\phi\|_{\infty}\int_0^T\intr (1-s)\int_0^T\intr |\Dv g^n|^2\Phi_{s,R}^n(g^n)dxd\vi dt\]
\[+\exp(sR)\int_0^T\intr \Phi_{s,R}^n(g^n)|\Dv g^n|^2\mathbf{1}_{g^n>R}.\]
Using the fact that
\[|\Phi_{s,R}^n(g^n)|\leq (1+m_n)^{s},\]
we obtain 
\[|\Dv g^n|^2\Phi_{s,R}^n(g^n)\leq \frac{|\Dv m^n|^2}{(1+m^n)^{2-s}}\leq \frac{|\Dv m^n|^2}{m^n} .\]
Furthermore,
\[\Phi_{s,R}^n(g^n)|\Dv g^n|^2\mathbf{1}_{g^n>R}\leq \exp(sR)\exp(-R)\frac{|\Dv m^n|^2}{m^n},\]
where in the last inequality we used that 
\[\Phi_{s,R}(t)=\exp(sR) \text{ for }t>R, \text{ and }\frac{1}{1+m^n}\mathbf{1}_{g^n>R}\leq \exp(-R).\]
Thus, from Proposition \ref{Log-EstimatesPropForm}, for some $C=C(\|m_0\|_1,\|m_0\log(m_0)\|_1,\|\log(1+m_0)\|_1,\sup\limits_{n}(\|b^n\|_2+\|m^n|b^n|^2\|_1))$ we have the bound
\[|(I)|\leq \Big( (1-s)\|\phi\|_{\infty}\]
\[+\exp(-(1-s)R)\Big)\int_0^T\intr \frac{|\Dv m^n|^2}{m^n}dxd\vi dt\leq C\Big( (1-s)\|\phi\|_{\infty}+\exp(-(1-s)R)\Big),\]
where in the last inequality is due to Proposition \ref{Log-EstimatesPropForm}. For the second term we work as follows
\[|(II)|\leq (1-s)\int_0^T\intr \Phi_{s,R}(g^n)\frac{|m^n||b^n|}{(1+m^n)}|\Dv g^n|dxd\vi dt\]
\[+\int_0^T\intr \Phi_{s,R}(g^n)\frac{|m^n||b^n|}{(1+m^n)}|\Dv g^n|\mathbf{1}_{m^n>R} dxd\vi dt.\]
For the first term above we use
\[\Phi_{s,R}(g^n)\frac{|m^n||b^n|}{(1+m^n)}|\Dv g^n|\leq \frac{|m^n||b^n|}{(1+m^n)^{2-s}}|\Dv m^n|\leq m_n|b^n|^2+\frac{|\Dv m^n|^2}{m^n},\]
while for the second integral
\[\Phi_{s,R}(g^n)\frac{|m^n||b^n|}{(1+m^n)}|\Dv g^n|\mathbf{1}_{g^n>R}\leq \exp(-(1-s)R)\Big(m^n|b^n|^2+\frac{|\Dv m^n|^2}{m^n} \Big),\]
hence
\[|(II)|\leq C\Big((1-s)+\exp(-(1-s)R)\Big).\]
Thus passing to the limit in \eqref{EqTestAgainPhin}, we obtain
\begin{equation}
    \label{FirstEstimateForMeasureVanish}
    \Big|\int_0^T \intr \phi'(t)\Psi_{s,R}(g)dxd\vi dt \Big|\leq C(1+\|\phi\|_{\infty})\Big((1-s)+e^{-(1-s )R} \Big).
\end{equation}
Now that we have obtained these bounds we obtain the result just as in \cite{diperna1988fokker}, section III. The only difference in the proof is the divergence term, which however causes no technical difficulties. We provide the details next. For $\epsilon>0$ let $\rho_{\epsilon}$ be a standard sequence of mollifiers. The equation satisfied by $g^{\epsilon}:=\rho_{\epsilon}\star g$ where $g$ solves
\[\pt g-\lapv g-\vDx g-\divv(\frac{m}{1+m}b)=|\Dv g|^2+\mu+\frac{m}{1+m}b\Dv g\]
reads
\begin{equation}
    \label{MollifiedForg}
    \pt \gep-\lapv \gep -\vDx \gep -\divv(\rho_{\epsilon}\star (\frac{m}{1+m}b))=\rho_{\epsilon}\star |\Dv g|^2+\rho_{\epsilon}\star (\frac{m}{1+m}b\Dv g)+\rho_{\epsilon}\star \mu+r_{\epsilon}.
\end{equation}
Testing against $\phi \Phi_{s,R}(\gep)$ in \eqref{MollifiedForg} yields,
\[-\int_0^T\intr \phi'(t)\Psi_{s,R}(\gep)dxd\vi dt\]
\[\geq \int_0^T\intr \phi(t)\Big[-|\Dv \gep|^2\Phi_{s,R}'(\gep)+\Phi_{s,R}'(\gep)\Dv \gep\rho_{\epsilon}\star (\frac{m}{1+m}b)+\rho_{\epsilon}\star |\Dv g|^2\Phi_{s,R}(\gep) +\]
\[\rho_{\epsilon}\star (\frac{m}{1+m}b\Dv g)\Phi_{s,R}(\gep)\Big] \phi(t)\Phi_{s,R}(\gep)\rho_{\epsilon}\star \mu dxd\vi dt-\|r_{\epsilon}\|_1\|\phi\|_{\infty}\|\Phi_{s,R}(\gep)\|_{\infty}.\]
We let $\epsilon\rightarrow 0$ and using that $\Phi_{s,R}\geq 1$ obtain
\[-\int_0^T\intr \phi'(t)\Psi_{s,R}(g)dxd\vi dt\]
\[\geq \int_0^T\intr \phi(t)\Big[|\Dv g|^2\Phi_{s,R}(g)-|\Dv g|^2\Phi_{s,R}'(g)\Big]\]
\[+\phi(t)\Big[\frac{m}{1+m}b\Dv g\Phi_{s,R}(g)-\Phi_{s,R}'(g)\Dv g\frac{m}{m+1}b \Big]dxd\vi +\int_0^T\intr \phi(t)d\mu\]
\[\geq \int_0^T\intr (1-s)\phi(t)|\Dv g|^2(g)\mathbf{1}_{g\leq R}+\phi(t)|\Dv g|^2\mathbf{1}_{g>R}\]
\[+\int_0^T\intr (1-s)\phi(t)\frac{m}{m+1}b\Dv g\Phi_{s,R}(g)\mathbf{1}_{g\leq R}+\phi(t)\frac{m}{m+1}b\Dv g\Phi_{s,R}(g)\mathbf{1}_{g>R},\]
where in the last equality we used that $\Phi_{s,R}\geq 1$. Next we bound the terms in the RHS 
\[\Big|\int_0^T\intr (1-s)\phi(t)|\Dv g|^2(g)\mathbf{1}_{g\leq R}+\phi(t)|\Dv g|^2\mathbf{1}_{g>R}\Big|\]
\[\leq (1-s)C\|\phi\|_{\infty}\|\Dv \sqrt{m}\|_2+\|\phi\|_{\infty}e^{-R}\|\Dv \sqrt{m}\|_2,\]
while for the rest of the terms 
\[\Big|\int_0^T\intr (1-s)\phi(t)\frac{m}{m+1}b\Dv g\Phi_{s,R}(g)\mathbf{1}_{g\leq R}+\phi(t)\frac{m}{m+1}b\Dv g\Phi_{s,R}(g)\mathbf{1}_{g>R} \Big|\]
\[\leq (1-s)\|\phi\|_{\infty}\Big(\|m|b|^2\|_1+\|\Dv \sqrt{m}\|_2 \Big)+\|\phi\|_{\infty}e^{-R(1-s)}\Big(\|m|b|^2\|_1+\|\Dv \sqrt{m}\|_2 \Big).\]
Hence combining the estimates above with estimate \eqref{FirstEstimateForMeasureVanish}, we obtain 
\[\int \phi d\mu\leq C((1-s)+e^{-R(1-s)})\]
letting $R\rightarrow \infty$ and then $s\uparrow 1$ yields $\int \phi d\mu\leq 0$, for all $\phi\geq 0$ and since $\mu\geq 0$ it follows that $\mu\equiv 0$. Finally, we show that $m\in C([0,T];L^1(\R^d\times \R^d))$. Let $\rho_{n}$ be a standard sequence of mollifiers (see section 1 for definition) and $m_n:=\rho_n\star m$. The functions $m_n$ satisfy 
\begin{equation}
    \label{MollifiedMForTimeCont}
    \pt m_n-\lapv m_n-\vDx m_n-\divv(\rho_n\star (mb))=r_n,\,\,
    m_n(0)=\rho_n\star m_0,
\end{equation}
where $r_n=K_n\star m$ and $K_n$ is given by
\[K_n:=n^{2d}\frac{\vi}{n}D_x \rho(\frac{x}{n},\frac{\vi}{n}),\]
and so $r_n\rightarrow 0$ strongly in $L^1([0,T]\times \R^d\times \R^d)$. From Lemma A.1 in \cite{degond1986global}, we have that $m_n\in C([0,T];L^2(\R^d\times \R^d))$. For any $S\in C_c^{\infty}(\R)$, the function $S(m_n)$ satisfies
\begin{equation*}
    \begin{cases}
    \pt S(m_n)-\lapv S(m_n)-\vDx S(m_n)-\divv(S'(m_n)\rho_n\star (mb))=-S''(m_n)|\Dv m_n|^2\\
    -S''(m_n)\Dv m_n\rho_n\star (mb)+S'(m_n)r_n,\\
    S(m_n)(0)=S(\rho_n\star m_0).
    \end{cases}
\end{equation*}
For $n,k\in \N$, we test against $S(m_n)-S(m_k)$ in the equation satisfied by their difference which yields for all $t\in [0,T]$
\[\intr |S(m_n)-S(m_k)|^2(t)dxd\vi+\int_0^t\intr |\Dv (S(m_n)-S(m_k))|^2dxd\vi dt \]
\[=-\boxed{\int_0^t\intr \Dv(S(m_n)-S(m_k))\Big(S'(m_n)\rho_n\star (mb)-S'(m_k)\rho_k\star (mb) \Big)dxd\vi dt}_1\]
\[-\boxed{\int_0^t\intr \Big( S(m_n)-S(m_k)\Big)\Big(S''(m_n)|\Dv m_n|^2-S''(m_k)|\Dv m_k|^2 \Big)dxd\vi dt}_2\]
\[-\boxed{\int_0^t\intr \Big( S(m_n)-S(m_k)\Big)\Big(S''(m_n)\Dv m_n\rho_n\star (mb)+S'(m_n)r_n-S''(m_k)\Dv m_k\rho_k\star (mb)}_3\]
\[\boxed{+S'(m_k)r_k \Big)dxd\vi dt }_3+\boxed{\intr |S(m_n)-S(m_k)|^2(0)dxd\vi}_4 .\]
As noted in \cite{porretta2015weak} (Remark 3.9) we have that 
\begin{equation}
    \label{BoundOnMollifier}
    |\rho_n \star (mb)|^2\leq [\rho_n\star (m|b|^2)]m_n.
\end{equation}
For the first boxed term the following hold 
\[\Dv (S(m_n))\rightarrow \Dv S(m) \text{ strongly in }L^2([0,T]\times \R^d\times \R^d) \text{ as }n\rightarrow \infty,\]
while from \eqref{BoundOnMollifier}, we obtain
\[|S'(m_n)\rho_n\star (mb)|^2\leq (S'(m_n))^2m_n[\rho_n\star (m|b|^2)]\leq C_S[\rho_n\star (m|b|^2)],\]
where $C_S:=\|(S'(x))^2x\|_{\infty}$. Since $[\rho_n\star (m|b|^2)]\rightarrow m|b|^2$ strongly in $L^1([0,T]\times \R^d\times \R^d)$ by Dominated Convergence Theorem we obtain 
\[S'(m_n)\rho_n\star (mb)\rightarrow S'(m)mb \text{ strongly in }L^2([0,T]\times \R^d\times \R^d),\]
therefore the first term can be bounded by a function $\omega(n,k)$ such that $\lim\limits_{n,k}\omega(n,k)=0$ independently of $t$. For the second term we note that 
\[S''(m_n)|\Dv m_n|^2\rightarrow S''(m)|\Dv m|^2 \text{ strongly in }L^1([0,T]\times \R^d\times \R^d),\]
while $S(m_n)\rightarrow S(m)$ strongly in $L^1([0,T]\times \R^d\times \R^d)$ with $\sup\limits_{n}\|S(m_n)\|_{\infty}<\infty$ therefore it can also be bounded like the first term. For the third term, from \eqref{BoundOnMollifier} we have
\[|S''(m_n)\Dv m_n\rho_n\star (mb)|\leq \frac{1}{2}|S''(m_n)||\Dv m_n|^2+|S''(m_n)m_n|[\rho_n\star (m|b|^2)]\]
and since the right hand side of the above inequality converges strongly in $L^1$ by Dominated Convergence we obtain that 
\[S''(m_n)\Dv m_n\rho_n\star (mb)\rightarrow S''(m)\Dv m\cdot mb \text{ strongly in }L^1([0,T]\times \R^d\times \R^d),\]
while $S'(m_n)r_n$ converges strongly to $0$ in $L^1([0,T]\times \R^d\times \R^d)$ just as in step 3, section III of \cite{diperna1988fokker}. Finally the fourth term clearly vanishes as $n,k\rightarrow \infty$. Thus taking the sup over $t$ we obtain 
\[\lim\limits_{t,n,k}\intr |S(m_n)-S(m_k)|^2(t)dxd\vi =0.\]
The above show that $S(m)\in C([0,T]; L^2(\R^d\times \R^d))$ for all $S\in C_c^{\infty}(\R^d\times \R^d)$ and so $T_k(m)\in C([0,T];L^2(\R^d\times \R^d))$ for all $k\in \N$ where $T_k$ is the truncation at $k$. To conclude, since for all $R>0$
\[\|m(t)-m(s)\|_{L^1(\R^d\times \R^d)}\leq \|m(t)-m(s)\|_{L^1(B_R)}+\|m(t)-m(s)\|_{L^1(B_R^c)}\]
and due to the bounds of Lemma \ref{Lemma|x|^2+|v|^2m}, we obtain that for some $C=C(R)>0$ and $C_1=C_1(m_0,b)>0$
\[\|m(t)-m(s)\|_{L^1(\R^d\times \R^d)}\leq C(R)\|T_k(m(t))-T_k(m(s))\|_2+2\sup\limits_{\theta \in [0,T]}\|m(\theta)-T_k(m(\theta))\|_1+\frac{C_1}{R^2}.\]
Furthermore by Proposition \ref{Log-EstimatesPropForm},
\[\|m(\theta)-T_k(m(\theta))\|_1=\int_{m(\theta)>k}|m|(\theta)dxd\vi \leq \frac{A(\|m_0\log(m_0)\|_1)}{\log(k)}, \]
where $A>0$ is the constant provided by Proposition \ref{Log-EstimatesPropForm}. Putting everything together we obtain 
\[\|m(t)-m(s)\|_1\leq C_R\|T_k(m(t))-T_k(m(s))\|_2+\frac{A}{\log(k)}\|m_0\log(m_0)\|_1+\frac{C_1}{R^2}.\]
Thus given an $\epsilon>0$, first we fix an $R>0$ such that $\frac{C_1}{R^2}\leq \frac{\epsilon}{3}$
and a $k\in \N$ such that 
\[\frac{A}{\log(k)}\|m_0\log(m_0)\|_1<\frac{\epsilon}{3},\]
then we find a $\delta >0$ such that \[|t-s|<\delta\implies C_R\|T_k(m(t))-T_k(m(s))\|_2<\frac{\epsilon}{3}\]
and so $m\in C([0,T]; L^1(\R^d\times \R^d))$.
\end{proof}
\subsection{Analysis of the Hamilton-Jacobi-Bellman equation}
In this section we will study the bounds for the HJB equation
\begin{equation}\label{SubsectionHJB}
\begin{cases}
-\pt u-\lapv u+\vDx u+H(\Dv u)=f(t,x,\vi) \text{ in }(0,T)\times \R^d\times \R^d,\\
u(T,x,\vi)=g(x,\vi) \text{ in }\R^d\times \R^d.
\end{cases}
\end{equation}
\begin{definition}\label{DefinitionOnCompactnessHJB}
Let $H:\R^d\rightarrow \R$ be a convex Lipschitz function such that $H\geq 0$, $f\in L^1\cap L^{\infty}([0,T]\times \R^d\times \R^d), f\geq 0, (|x|^2+|\vi|^2)f\in L^1([0,T]\times \R^d\times \R^d) \,\,g\in L^1\cap L^{\infty}(\R^d\times \R^d), g\geq 0, (|x|^2+|\vi|^2)g\in L^1(\times \R^d\times \R^d)$ and $u\in C([0,T];L^2(\R^d\times \R^d))\cap L^1(\R^d\times \R^d)$ with $\Dv u\in L^2([0,T]\times \R^d\times \R^d),u\geq 0$. We say that $(u,H,f,g)$ is a weak solution of \eqref{SubsectionHJB}, if the equation is satisfied in the distributional sense.

\end{definition}
Our starting point is the following compactness theorem found in the Appendix of  \cite{diperna1988fokker}.
\begin{theorem}[Appendix of P.-L. Lions, DiPerna \cite{diperna1988fokker}]\label{L1-comp}
Assume that $u^n,f^n\in L^1([0,T]\times \R^d\times \R^d),g^n\in L^1(\R^d\times \R^d)$ satisfy in the distributional sense
\begin{equation*}
\pt u_n-\lapv u_n+\vDx u_n=f_n,\,\, u_n(0)=g^n.
\end{equation*}
If $g^n,f_n$ are uniformly bounded in $L^1$ with 
 \begin{equation}
 \label{Vanishinggn}
 \lim\limits_{R\rightarrow \infty}\sup\limits_{n\in \N}\int_0^T\int_{|(x,\vi)|\geq R}|f^n|dxd\vi dt=0
 \end{equation}
and
 \begin{equation}
 \label{Vanishingg0n}
 \lim\limits_{R\rightarrow \infty}\sup\limits_{n\in \N}\int_{|(x,\vi)|\geq R}|g_0^n|dxd\vi =0,
 \end{equation}
  then the sequence $\{u_n\}_{n\in \N}$ is compact in $L^1((0,T)\times \R^d\times \R^d)$. 
\end{theorem}
\begin{theorem}\label{VanishingHJBComp}
Let $f^n\in L^1([0,T]\times \R^d\times \R^d), g^n\in L^1(\R^d\times \R^d)$ be non-negative, uniformly integrable sequences and $H^n:\R^d\rightarrow \R$ Lipschitz convex Hamiltonians. Assume that $\{(u^n,H^n,f^n,g^n)\}_{n\in \N}$ are weak solution to \eqref{SubsectionHJB} according to definition \ref{DefinitionOnCompactnessHJB}. Then, the sequence $\{u^n\}$ is compact in $L^1((0,T)\times \R^d\times \R^d)$ and 
\[\sup\limits_{n\in \N}\Big(\sup\limits_{t\in [0,T]}\|u^n(t)\|_1+\|H^n(\Dv u^n)\|_1 \Big)<\infty,\]
\[\lim\limits_{R\rightarrow \infty}\sup\limits_{n}\Big(\sup\limits_{t\in [0,T]}\int_{B(0,R)^c}|\un|(t)dxd\vi +\int_{B(0,R)^c}H^n(\Dv \un)dxd\vi dt \Big)=0.\]
\end{theorem}
\begin{proof}
By the same arguments as in Lemma \ref{Lemma|x|^2+|v|^2m}, we can justify testing against $1$ in the HJB equation to obtain the uniform $L^1$ estimates on $\un,H^n(\Dv \un)$. To show compactness in $L^1$ we work as follows. Let $L:=-\pt -\lapv+ \vDx$ and since $H^n\geq 0,f^n\geq 0,g^n\geq 0$ the functions $\un$ are non-negative and satisfy 
\begin{equation*}
    L\un\leq f^n \text{ in }(0,T)\times \R^d\times \R^d,\,\,
    \un(T)=g^n \text{ in }\R^d\times \R^d.
\end{equation*}
For each $n\in \N$, let $w^n$ be the solution of 
\begin{equation*}
    Lw^n=f^n \text{ in }(0,T)\times \R^d\times \R^d,\,\,
    w^n(T)=g^n \text{ in }\R^d\times \R^d.
\end{equation*}
Since $L(w^n-\un)\geq 0$ and $w^n(T)=\un(T)$ we have that 
\begin{equation}\label{IntermediateSolutions}
    0\leq \un\leq w^n.
\end{equation}
Since $f^n,g^n$ are uniformly integrable, by Theorem \ref{L1-comp} the set $\{w^n\}_{n\in \N}$ is compact in $L^1$ and so in particular uniformly integrable and from \eqref{IntermediateSolutions} we see that $\{\un\}_{n\in \N}$ are also uniformly integrable. For $R>0$, let $\phi_R:\R^d\times \R^d\rightarrow [0,1]$ be cutoff functions defined just as in Lemma \ref{Lemma|x|^2+|v|^2m}. Testing against $\phi_R$ in
\begin{equation*}
    L\un+H(\Dv \un)=f^n,\,\,\un(T)=g^n
\end{equation*}
yields for some dimensional constant $C>0$
\[\intr \un(t)\phi_Rdxd\vi +\int_0^t\intr H^n(\Dv \un)\phi_Rdxd\vi dt\leq\]
\[\frac{C}{R}\|\un\|_1+\int_0^t\intr f^n\phi_Rdxd\vi dt+\int_0^t\intr g^n\phi_R+C\int_{R<|(x,\vi)|<2R}\un dxd\vi \]
and since the sequence $\{\un\}_{n\in \N}$ is uniformly integrable we see that the terms on the right vanish uniformly in $n$ as $R\uparrow\infty$. Finally with the estimate 
\[\lim\limits_{R\rightarrow \infty}\sup\limits_{n\in \N}\int_{R<|(x,\vi)|}H^n(\Dv \un)dxd\vi dt=0\]
the compactness for $\un$ in $L^1$ follows immediately by Theorem \ref{L1-comp} with $\tilde{f}^n=f^n-H^n(\Dv \un)$.
\end{proof}

\begin{theorem}\label{GradientAndL2estimatesU}
Let $(u,H,f,g)$ be a weak solution of \eqref{SubsectionHJB}, according to Definition \ref{DefinitionOnCompactnessHJB}. Then, there exists a constant $C=C(d,T)>0$, such that
\begin{equation}
\label{Duepsilon^2}
\sup\limits_{t\in [0,T]}\|u(t)\|_2+\|u H(\Dv u)\|_1+\|\Dv u\|_2\leq C\Big( \|f\|_{\infty}\|f\|_1+\|g\|_1\|g\|_{\infty}\Big).
\end{equation}
\end{theorem}
\begin{proof}
The result follows by testing against $u$ in \eqref{SubsectionHJB} and applying Gr\"onwall.
\end{proof}
\begin{proposition}\label{ConvergenceOfGradientDu}
Let $\{(u^n,H^n,f^n,g^n)\}_{n\in\N}$, be weak solutions of \eqref{SubsectionHJB}, according to Definition \ref{DefinitionOnCompactnessHJB}, such that 
\[\|f^n\|_1+\|g^n\|_1\leq C\text{ for all }n\in\N,\]
and 
\[u^n\rightarrow u \text{ strongly in }L^1([0,T]\times \R^d\times \R^d).\]
Then, the limit $u$ belongs to $L^2([0,T]\times \R^d; H^1(\R_v^d))$ and 
\[\Dv u^n\rightarrow \Dv \text{ in }L_{loc}^q([0,T]\times \R^d\times \R^d),\]
for all $q<2$, up to a subsequence almost everywhere. 
\end{proposition}
\begin{proof}
The equation for $u^n-u^m$ is
\[-\pt(u^n-u^m)-\lapv (u^n-u^m)+\vDx (u^n-u^m)=f^n-f^m,\]
\[(u^n-u^m)(T)=g^n-g^m.\]
For $\epsilon>0$, we define
\begin{equation*}
\phi(s):=\begin{cases}
s, \text{ for }s\in [-\epsilon,\epsilon],\\
-\epsilon, \text{ for }s\leq -\epsilon,\\
\epsilon, \text{ for }s\geq \epsilon,
\end{cases}
\end{equation*}
and $\Phi(t):=\int_0^t \phi(s)ds\geq 0$. We test against $\phi(u^n-u^m)$ in the equation for the difference, which yields
\[\intr \Phi(u^n-u^m)(t)dxd\vi +\int_0^T\intr \phi'(u^n-u^m) |\Dv (u^n-u^m)|^2dxd\vi\]
\[ \leq \intr \Phi(u^n-u^m)(T)dxd\vi + \int_0^T\intr \phi(u^n-u^m)(f^n-f^m)dxd\vi dt\]
\[\leq C\epsilon \|g^n-g^m\|_1+\epsilon\|f^n-f^m\|_1\leq C\epsilon.\]
Therefore,
\[\int_{|u^n-u^m|\leq \epsilon}|\Dv(u^n-u^m)|^2dxd\vi dt\leq C\epsilon.\]
Thus, fixing a radius $R>0$ and a $q<2$ we obtain
\[\int_{B(0,R)} |\Dv (u^n-u^m)|^qdxd\vi\leq \int_{B(0,R)\cap \{|u^n-u^m|\leq \epsilon\}}|\Dv(u^n-u^m)|^qdxd\vi dt\]
\[+\int_{B(0,R)\cap |u^n-u^m|>\epsilon}|\Dv(u^n-u^m)|^qdxd\vi dt \leq CR^d\epsilon+CR^d|\{|u^n-u^m|>\epsilon\}|^{\theta}\]
for some $\theta=\theta(q)\in (0,1)$. Since $u^n$ converges in $L^1$, we have that $\lim\limits_{n,m\rightarrow \infty}|\{|u^n-u^m|>\epsilon\}|=0$
and so $\Dv u^n\rightarrow \Dv u \text{ in }L^q([0,T]\times B(0,R)) \text{ for all }R>0.$
\end{proof}
\begin{proposition}\label{ConvergenceOfH}
Assume that $\{(u^n,H^n,f^n,g^n)\}_{n\in \N}$ are weak solutions to \eqref{SubsectionHJB} according to Definition \ref{DefinitionOnCompactnessHJB}, such that $\{g^n \}_{n\in \N}\subset L^1(\R^d\times \R^d)$ is uniformly integrable, $\{f^n\}_{n\in \N}$ and $\{g^n\}_n$ are bounded subsets of their respective $L^{\infty}$ spaces, and for some $u,f$, $u^n\rightarrow u, f^n\rightarrow f, f^n\rightarrow f$, in $L^1([0,T]\times \R^d\times \R^d)$ and almost everywhere. Then, up to a subsequence, for each $\tau \in [0,T)$, we have that 
 \[H^n(\Dv u^n)\rightarrow H(\Dv u)\text{ in }L^1([0,\tau]\times \R^d\times \R^d)\]
and,
\[\Dv u^n\rightarrow \Dv u \text{ in }L^2([0,\tau]\times \R^d\times \R^d).\] 
\end{proposition}
\begin{proof}
From Proposition \ref{ConvergenceOfGradientDu} by choosing a subsequence if necessary we can assume that $H^n(\Dv \un)\rightarrow H(\Dv u)$ almost everywhere, furthermore since $\sup\limits_{n}\|f^n\|_{\infty}+\|g^n\|_{\infty}<\infty$, for some $C>0$ we have that $\|\un\|_{\infty}\leq C$ for all $n\in \N$. Denote by $L:=-\pt -\lapv +\vDx$ and test against $(T-t)e^{\lambda(\un-u^k)}$ in the equation
\[L(\un-u^k)+[H^n(\Dv \un)-H^k(\Dv u^k)]=f^n-f^k.\]
Which yields,
\[\intr T\frac{1}{\lambda}(e^{\lambda(\un-u^k)}-1)(0)dxd\vi -\int_0^T\intr \frac{1}{\lambda}(e^{\lambda(\un-u^k)}-1)(s)dxd\vi ds\]
\[+\int_0^T\intr (T-s) e^{\lambda(\un-u^k)}|\Dv (\un-u^k)|^2+(T-s) e^{\lambda(\un-u^k)}\Big(H^n(\Dv \un)-H^k(\Dv u^k) \Big)dxd\vi ds\]
\[=\int_0^T\intr e^{\lambda(\un-u^k)}\Big(f^n-f^k \Big)dxd\vi ds.  \]
Next using the strong convergence of $\un,f^n$ and that $\un$ is uniformly bounded in $L^{\infty},$ we obtain that for some function $\omega(n,k)$ such that $\lim\limits_{n,k\rightarrow \infty}\omega(n,k)=0$
\[\int_0^T\intr (T-s)\lambda e^{\lambda(\un-u^k)}|\Dv (\un-u^k)|^2dxd\vi ds\]
\[+\int_0^T\intr (T-s)e^{\lambda(\un-u^k)}\Big(H^n(\Dv \un)-H^k(\Dv u^k) \Big)dxd\vi ds\leq \omega(n,k)\]
If $n>k$ we have that $H^k\leq H^n$, hence by the convexity of $H$
\[\int_0^T\intr (T-s) \lambda e^{\lambda(\un-u^k)}|\Dv (\un-u^k)|^2dxd\vi ds\]
\[+\int_0^T\intr(T-s) e^{\lambda(\un-u^k)}\Big(H^n(\Dv \un)-H^n(\Dv u^k) \Big)dxd\vi ds\leq \omega(n,k)\]
\[\implies \int_0^T\intr (T-s)\lambda e^{\lambda(\un-u^k)}|\Dv (\un-u^k)|^2dxd\vi ds\]
\[+\int_0^T\intr(T-s) e^{\lambda(\un-u^k)}H_p^n(\Dv u^k)\Dv(\un-u^k)dxd\vi ds\leq \omega(n,k).\]
Letting $n\rightarrow \infty$ and using that $\Dv \un\rightarrow \Dv u$ almost everywhere and weakly in $L^2$, while $\un\rightarrow u$ strongly in $L^1$ with $\|\un\|_{\infty}\leq C$ and $|H_p^n(\Dv u^k)|\leq |H_p|(\Dv u^k)$ thus $H_p^n(\Dv u^k)\rightarrow H_p(\Dv u^k)$ strongly in $L^2$, yields
\[\int_0^T\intr (T-s)\lambda e^{\lambda(u-u^k)}|\Dv (u-u^k)|^2dxd\vi ds\]
\[+\int_0^T\intr(T-s) e^{\lambda(u-u^k)}H_p(\Dv u^k)\Dv(u-u^k)dxd\vi ds\leq \omega(k).\]
From \ref{Assumption2onH}, there exists a constant $C>0$ such that 
\[H_p(\Dv u^k)\Dv (u-u^k)=-(H_p(\Dv u)-H_p(\Dv u^k))\cdot \Dv (u-u^k)+H_p(\Dv u)\Dv (u-u^k)\]
\[\geq -C |\Dv (u-u^k)|^2+H_p(\Dv u)\Dv (u-u^k)\]

\[\implies \int_0^T\intr (T-s) e^{\lambda(u-u^k)}(\lambda-C)|\Dv (u-u^k)|^2dxd\vi ds\]
\[+\int_0^T\intr(T-s) e^{\lambda(u-u^k)}H_p(\Dv u)\Dv(u-u^k)dxd\vi ds\leq \omega(k)\]
and again by the weak convergence of $\Dv (u-u^k)$ in $L^2$ and the strong convergence of $u^k$ to $u$ in $L^1$ with uniform bounds we obtain
\[\int_0^T\intr (T-s) e^{\lambda(u-u^k)}(\lambda-C)|\Dv (u-u^k)|^2dxd\vi ds\leq \omega(k).\]
Finally, the result follows since by choosing $\lambda>C$ and using that $\|u-u^k\|_{\infty}\leq C$ we obtain that for some constant $c_0>0$ depending only on $H$ 
\[c_0\int_0^T\intr (T-s) |\Dv (u-u^k)|^2dxd\vi ds\leq \omega(k).\]
\end{proof}

\begin{theorem}\label{L1-time-convergenceForU}
Assume that $\{(u^n,H^n,f^n,g^n)\}_{n\in \N}$ are weak solutions to \eqref{SubsectionHJB} according to Definition \ref{DefinitionOnCompactnessHJB}, such that $f^n\rightarrow f$ in $L^1$, $g^n\rightarrow g$, weakly in $L^1$, $u^n\rightarrow u$ in $L^1$ and $\Dv u^n\rightarrow \Dv u$ almost everywhere and $H^n(\Dv u^n)\rightarrow H(\Dv u)$ in $L_{loc}^1((0,T];L^1(\R^d\times \R^d))$, where $H(\Dv u)\in L^1([0,T]\times \R^d\times \R^d)$. Then, we have that $u\in C((0,T];L^1((\R^d\times \R^d))$.
\end{theorem}
\begin{proof}
The result follows by the fact that $Lu\in L^1$, where $L:=-\pt -\lapv +\vDx$, see for example \cite{diperna1988fokker}.
\end{proof}
\subsection{Existence and uniqueness for the quadratic case}
In this subsection, we will establish the existence and uniqueness of renormalized solutions for the MFG system.
\begin{theorem}\label{ExistenceAndUniquenessQuadratic}
Assume that $H:\R^d\rightarrow \R,\,\, F:[0,T]\times \R^d\times \R^d\times \R\rightarrow \R,\,\, m_0:\R^d\times \R^d\rightarrow \R$ and $G:\R^d\times \R^d\times \R\rightarrow \R$ satisfy \ref{Assumption2onH},\ref{F_LAssum},\ref{AssumptionInitial} and \ref{G_LAssum}. Then, there exists a unique renormalized solution $(m,u)$ of system \eqref{MFGequation}, according to Definition \ref{DefinitionOfRenormalizedToMFGSystem}.
\end{theorem}
\begin{proof}
The proof is divided in two steps. First we show the result for $F,G$ bounded in their respective $L^{\infty}-$spaces and let the Hamiltonians $H^{\epsilon}$ vary. While in the second case we show the result for a fixed quadratic Hamiltonian $H$ while letting $F^n,G^n$ vary. The reason for this approach is so that we can always have bounds on $\Dv \un$ in $L^2$. Indeed in the first case the bounds follow by Theorem \ref{GradientAndL2estimatesU} and are due to the $\lapv $ term while in the second case the bounds are a result of Theorem \ref{VanishingHJBComp} and are due to $\|H(\Dv \un)\|_1\leq C$.
\paragraph{First Case:}
For $H^{\epsilon}$, as defined in \eqref{Hepsilon}, we consider the solutions $(\me,\ue,m_0)$ provided by Theorem \ref{ExistenceTheorem}. From Proposition \ref{EnergyEstimateProp} above, we have that for some $C>0$ independent of $\epsilon$ 
\begin{equation}
\label{FUNDAMENTALFOREXISTENCE}
\|\me |H_p^{\epsilon}(\Dv \ue)|^2\|_1\leq C, \text{ for all }\epsilon>0,
\end{equation}
furthermore, by Theorem \ref{GradientAndL2estimatesU} and our assumptions on $H^{\epsilon}$ we have that
\[\|H_p^{\epsilon}(\Dv \ue)\|_2\leq C, \text{ for all }\epsilon>0.\]
Therefore, from Theorem \ref{CompactnesL^1Form}, we may extract a subsequence $m^n$, which is convergent in $L^1([0,T]\times\R^d\times \R^d)$ and almost everywhere to some $m$. From Remark \ref{RemarkUniformIntegralibity}, we have that the sequence $\{F(t,x,\vi,m^n)\}_{n\in \N}$ is uniformly integrable, indeed in the case $f_L:=\sup\limits_{m\in [0,L]}F(t,x,\vi,m)\in L^1$ the claim holds just as \cite{porretta2015weak}, while in the case $f_L:=\sup\limits_{m\in [0,L]}F(t,x,\vi,m)/m\in L^{\infty}$ since
\[0\leq F(t,x,\vi,m^n)\leq f_L(t,x,\vi) m^n+\frac{m^n}{L}F(t,x,\vi,m^n)\]
the result follows due to uniform bound on $\|F(t,x,\vi,m^n)m^n\|_1$ from Proposition \ref{EnergyEstimateProp} and the convergence of $m^n$ in $L^1$.
Since $m^n\rightarrow m$ almost everywhere, we obtain 
\[F(\cdot,m^n(\cdot))\rightarrow F(\cdot ,m(\cdot)) \text{ strongly in }L^1([0,T]\times \R^d\times \R^d).\]
By choosing a further subsequence if necessary, Theorem \ref{VanishingHJBComp}, Lemma \ref{Lemma|x|^2+|v|^2m} and Proposition \ref{ConvergenceOfGradientDu}, yield a $u\in C([0,T];L^1( \R^d\times \R^d))\cap L^2([0,T]\times \R^d;H^1(\R_{\vi}^d))$, such that 
\[u^n\rightarrow u \text{ almost everywhere and strongly in }L^1([0,T]\times \R^d\times \R^d)\]
\[\Dv u^n\rightarrow \Dv u \text{ almost everywhere and in }L_{loc}^1([0,T]\times \R^d\times \R^d).\]
Furthermore, again by taking subsequences if needed, by Proposition \ref{ConvergenceOfH} we have that for each $\tau \in [0,T)$,
\[H^{\epsilon_n}(\Dv u^n)\rightarrow H(\Dv u) \text{ in }L^1([0,\tau]\times \R^d\times \R^d)\]
and for each $k\in \N$,
\[\Dv (u^n\wedge k)\rightarrow \Dv (u\wedge k) \text{ in }L^2([0,\tau]\times \R^d\times \R^d).\] 
By inequality \eqref{FUNDAMENTALFOREXISTENCE} and the fact that $H_p^{\epsilon_n}(\Dv u^n)\rightarrow H_p^{\epsilon_n}(\Dv u)$ almost everywhere, Proposition \ref{WeakLimitFormQuadratic} implies that
\[m^n\rightarrow m \text{ in }C([0,T];\mathcal{P}(\R^d\times \R^d)),\]
and by Theorem \ref{RenormalizedSolutionForm}, $m$ is a renormalized solution of
\begin{equation}
\label{PartOfProofEquationForLimitOfm}
\pt m-\lapv m-\vDx m-\divv(mH_p(\Dv u))=0 \text{ in }(0,T]\times \R^d\times \R^d,\\
m(0)=m_0 \text{ in }\R^d\times \R^d.
\end{equation}
It remains to show the convergence of the terminal data in the HJB equation. This follows exactly as in \cite{porretta2015weak}. Thus, we have that $m^n(T)\rightarrow m(T) \text{ in }L^1(\R^d\times \R^d)$
which from Remark \ref{RemarkUniformIntegralibity} implies that $G(\cdot,m^n(T,\cdot))\rightarrow G(\cdot,m(T,\cdot)) \text{ in }L^1(\R^d\times \R^d).$ Thus, the limit $u$ is also a renormalized solution.
\paragraph{Second Case:} Next, given $F,G$ that satisfy \ref{F_LAssum} and \ref{G_LAssum} respectively, consider $F^n:=F\wedge n, G^n:=G\wedge n$ for $n\in \N$. The functions $F^n,G^n$ clearly also satisfy \ref{F_LAssum} and \ref{G_LAssum} respectively. Let $(\un,m^n)$ be the solutions provided for the data $(H,F^n,G^n)$ by the first case. The rest of the proof follows exactly the first case only now we use Theorem \ref{ConvergenceOfGradientDuThroughPhia} to obtain the convergence of $\Dv T_k(\un)$.\\
Finally, we address the issue of uniqueness whose proof follows the same exact arguments as \cite{porretta2015weak} once we establish that $m(t,x,\vi)>0$ almost everywhere. But this will follow from assumption \eqref{AssumptionInitial} and in particular $\log(m_0)\in L_{loc}^1(\R^d\times \R^d)$. Indeed, let $R>0$ and define $\phi_R:\R^d\times \R^d\rightarrow [0,1]$ such that
\begin{equation*}
    \phi_R(x,\vi):=
    \begin{cases}
    1 \text{ if }|(x,\vi)|\leq R,\\
    0 \text{ if }|(x,\vi)|\geq R+1.
    \end{cases}
\end{equation*}
Then since 
\[\intr \log(m(t))\phi_R^2dxd\vi \leq \intr m(t)\phi_R^2dxd\vi \leq 1,\]
it is enough to bound $\intr \log(m(t))\phi_Rdxd\vi$ from bellow, since that would imply $m(t,x,\vi)>0$ almost everywhere. To show the lower bound we test the equation satisfied by $m$ with $\phi_R^2\frac{1}{m}$ (technically we would need to fix a $\delta >0$ and test against $\phi_R^2\frac{1}{m+\delta}$ and let $\delta\rightarrow 0$ but we skip the approximation for simplicity). This yields
\[ \intr \log(m(t))\phi_R^2 dxd\vi +\int_0^t\intr -\frac{|\Dv m|^2}{m^2}\phi_R^2+2\frac{\Dv m}{m}\phi_R\Dv \phi_R-\frac{\Dv m}{m^2}mH_p\phi_R^2 \]
\[+2\phi_R\Dv \phi_R H_p dxd\vi dt=\intr \log(m_0)\phi_R^2dxd\vi. \]
Next we use the following inequalities
\[2\Big|\frac{\Dv m}{m}\phi_R\Dv \phi_R\Big|\leq \frac{1}{4}\frac{|\Dv m|^2}{m^2}\phi_R^2+4|\Dv \phi_R|^2\]
\[\Big| \frac{\Dv m}{m^2}mH_p\phi_R^2\Big|\leq \frac{1}{4}|\frac{\Dv m}{m}|^2\phi_R^2+|H_p|^2\phi_R^2\]
\[\Big| 2\phi_R\Dv \phi_R H_p\Big|\leq |H_p|^2+2|\phi_R|^2|\Dv \phi_R|^2\]
and thus combining everything we obtain that for some constant $C=C(R,d)>0$
\[\intr \log(m(t))\phi_R^2  \geq \intr C(R,d)+\|\log(m_0)\phi_R^2\|_1-\|H_p(\Dv u)\|_2\]
which proves the claim.
\end{proof}
\subsection{Further regularity for quadratic Hamiltonian} 
\begin{theorem}
Let $(H,F,G,m_0)$ be as in Theorem \ref{ExistenceAndUniquenessQuadratic} with $F=F(m)$, $G=G(m)$ and $m_0$ also satisfying $\|D^2m_0\|_{\infty}\in L^{\infty}(\R^d\times \R^d)$. Then, there exists a constant $C(F,G,H,T,m_0)$, such that 
\[\intr G'(m(T,x,\vi))|Dm(T,x,\vi)|^2dxd\vi +\int_0^T\intr F'(m(t,x,\vi))|Dm(t,x,\vi)|^2\]
\[+m\sum\limits_{k=1}^{2d}m\Dv u_{k}H_{pp}(\Dv u)\Dv u_k dxd\vi \leq C.\]
\end{theorem}
\begin{proof}
The proof is almost identical to the one in the case of Lipschitz Hamiltonian. The only difference is now instead of using the HJB equation we estimate
\[\intr \delta^hm_0\delta^h u(0)dxd\vi =\intr \frac{m_0(x+h,\vi)-2m_0(x,\vi)+m_0(x-h,\vi)}{h^2 }u(0,x,\vi)dx\vi\]
\[\leq \|D^2m_0\|_{L^{\infty}(\R^d\times \R^d)}\|u(0,\cdot,\cdot)\|_{1},\]
and conclude due to the estimate in Proposition \ref{EnergyEstimate}.
\end{proof}
\section{Appendix}
\subsection{Technical results}
In this sub-section we show some important properties about the convergence of $\un$ where $\un$ solves
\begin{equation}\label{AppendixEquation}
    \begin{cases}
    L\un+H(\Dv \un)=f^n \text{ in }(0,T]\times \R^d\times \R^d,\\
    \un(0)=g^n \text{ in }\R^d\times \R^d,
    \end{cases}
\end{equation}
for $L:=\pt -\lapv +\vDx$ and $f^n,g^n$ strongly convergent sequences in their respective $L^1$-spaces. We show an analogue of the convergence results in \cite{porretta1999existence} from which all our techniques are motivated and parallel to. In particular we show that if $\un$ solves \eqref{AppendixEquation} and are strongly convergent in $L^1$ to some function $u$, then $\Dv T_k(\un)\rightarrow \Dv T_k(u)$ strongly in $L^2([0,T]\times \R^d\times \R^d)$, where $T_k$ is the truncation at $k$, namely, for $k\in \N$ $T_k(x)=s$ for $|s|\leq k$ and $T_K(s)=sign(s)k$ otherwise. A crucial technical step in \cite{porretta1999existence} is the following transformation which allows the authors to deal with the degenerate $\pt$ direction. Given a function $u$, for $\nu>0$ define  
\[\pt (u)_{\nu}=\nu(T_k(u)-(u)_{\nu}).\] 
This transformation enjoys many nice properties such as $(u)_{\nu}\rightarrow u$ and $D((u)_{\nu})\rightarrow Du$ as $\nu\rightarrow \infty$ in appropriate spaces. In our setup the above transformation does not seem to work due to the extra degenerate operator $\vDx $. In order to deal with this, we consider a slightly different transformation. Fix $\alpha>0$ and consider the solution of
\[L\Phia=\alpha(T_k(u)-\Phia).\]
We will show that under the condition $u\in L^1$ the transformation $\Phia$ converges to $T_k(u)$ in $L^1$, however, we cannot show in general, even if $\Dv u\in L^2$, that $\Dv \Phia\rightarrow \Dv T_k(u)$ strongly in $L^2$, with no assumptions on $D_x u$. However the fact that $L\un +H(\Dv \un)=f^n$ and $\un\rightarrow u$ strongly in $L^1$, is enough to show the strong convergence of $\Dv \Phia$. With this, we can follow the rest of the argument found in \cite{porretta1999existence}. 
\begin{lemma}\label{PhiaConvergence}
Let $u\in L^1\cap L^{\infty}([0,T]\times \R^d\times \R^d)\cap C([0,T]; L^1(\R^d\times \R^d))$ and $\alpha>0$. Then, there exists a unique function $\Phia\in L^2([0,T]\times \R^d\times \R^d)$ with $\Dv \Phia \in L^2([0,T]\times \R^d\times \R^d)$ which solves
\begin{equation}
    \begin{cases}\label{L-RegularizationInLemma}
    \pt \Phia-\lapv \Phia +\vDx \Phia =\alpha(u-\Phia) \text{ in }(0,T)\times \R^d\times \R^d,\\
    \Phia(0,x,\vi)=u(0,x,\vi) \text{ in }\R^d\times \R^d.
    \end{cases}
\end{equation}
Furthermore, the functions $\Phia$ have the following properties
\begin{enumerate}
    \item\label{Positive} $u\geq 0\implies \Phia\geq 0$ almost everywhere,
    \item\label{UpperBounds} $\|\Phia \|_{\infty}\leq \|u\|_{\infty}$,\\
    \item\label{L2ConvergenceWholetime} $\lim\limits_{\alpha\rightarrow \infty}\|\Phia-u\|_2=0$\\
    \item $\|\Phia\|_1\leq \|u\|_1+\frac{1}{\alpha}\|u_0\|_1$\\
\end{enumerate}
\end{lemma}
\begin{proof}
First we assume that $u\in C^{\infty}([0,T];C_c^{\infty}(\R^d\times \R^d))$. Let $\Gamma$ denote the fundamental solution of $L=\pt -\lapv +\vDx $. Then, it is easy to check that the solution of equation \eqref{L-RegularizationInLemma} is given by 
\[\Phia(t,x,\vi)=\int_0^t\intr \alpha e^{-\alpha(t-s)}\Gamma(t-s,x,\vi,y,w)u(s,y,w)dydw ds\]
\[+\intr \alpha e^{-\alpha t}\Gamma(t,x,\vi,y,w)u(0,y,w)dydw ,\]
see for example \cite{lunardi1997schauder}. Furthermore, the solution $\Phia$ is also $C^{\infty}$ since $L$ is hypoelliptic. Let $f:=L(u)\in C^{\infty}([0,T]\times C_c( \R^d\times \R^d))$. In the equation 
\[L(u-\Phia)=-\alpha(u-\Phia)+f,\,\,\,(u-\Phia)(0)=0,\]
we test against $(u-\Phia)$, which yields
\[\frac{d}{dt}\frac{1}{2}\intr |u-\Phia|^2dxd\vi  +\intr |\Dv (u-\Phia)|^2dxd\vi\]
\[=-\alpha \intr |u-\Phia|^2dxd\vi +\intr f(u-\Phia)dxd\vi  \leq \frac{1}{4\alpha }\intr |f|^2dxd\vi .\]
Hence, we obtain that 
\[\sup\limits_{t\in [0,T]}\|u(t)-\Phia(t)\|_2+\|\Dv (u-\Phia)\|_2\leq \frac{C}{\alpha}\]
where $C=C(T,f)>0$. Furthermore, by testing against $p|u-\Phia|^{p-2}(u-\Phia)$ for $p>1$ yields
\[\frac{d}{dt}\int |u-\Phia|^pdxd\vi +\intr |\Dv(u-\Phia)|^2|u-\Phia|^{p-2}p(p-1)dxd\vi\]
\[\leq -\alpha p\intr |u-\Phia|^p +p\intr |f||u-\Phia|^{p-1}dxd\vi \leq \frac{p}{4a}\intr |f|^pdxd\vi, \]
where $1/p+1/q=1$. Letting $p\rightarrow 1$ yields
\[\sup\limits_{t\in [0,1]}\|u-\Phia\|_1\leq \frac{C}{\alpha}\|f\|_1,\]
where $C=C(T,f)>0$. The first two claims now follow easily by the Maximum Principle. For the general case we work as follows. Testing against $p|\Phia|^{p-2}\Phia$ in \eqref{L-RegularizationInLemma} for $p>1$ and letting $p\rightarrow 1$ just as above we obtain 
\[-\intr |u_0|dxd\vi \leq \frac{\alpha}{2}\int_0^T\intr |u|dxd\vi dt-\frac{\alpha}{2}\int_0^T\intr |\Phia|dxd\vi dt.\]
Hence,
\[\|\Phia\|_1\leq \|u\|_1+\frac{2}{\alpha}\|u_0\|_1,\]
and so by linearity of the map $(u,u_0)\rightarrow \Phia$ and the fact that $|u|\leq k\implies |\Phia|\leq k$ the result holds in the general case.
\end{proof}
\begin{theorem}\label{ConvergenceOfGradientDuThroughPhia}
Let $H:\R^d\rightarrow \R$ be a Hamiltonian satisfying \ref{Assumption2onH}. Assume that $\{f^n\}_{n\in \N}\subset L^1\cap L^{\infty}([0,T]\times \R^d\times \R^d), \{g^n\}_{n\in \N}\subset L^1\cap L^{\infty}([0,T]\times \R^d\times \R^d)$ such that $f^n\rightarrow f$ and $g^n\rightarrow g$ strongly in the respective $L^1$ spaces (the limits need not be in $L^{\infty}$). Let $\un \in L^1\cap L^2([0,T]\times \R^d\times \R^d)$ with $\Dv \un\in L^2([0,T]\times \R^d\times \R^d)$ solve
\begin{equation}
    \begin{cases}\label{HJBInPhiaConvergence}
    \pt \un -\lapv \un +\vDx \un +H(\Dv \un )=f^n, \text{ in }(0,T)\times \R^d\times \R^d,\\
    \un(0,x,\vi)=g^n(x,\vi) \text{ in }\R^d\times \R^d.
    \end{cases}
\end{equation}
Finally, assume that $\un\rightarrow u$ strongly in $L^1$ and that $\Dv \un \rightarrow \Dv u$ almost everywhere. Then, the limit $u$ is a renormalized solution of
\begin{equation*}
    \begin{cases}
    \pt u-\lapv u+\vDx u+H(\Dv u)=f(t,x,\vi) \text{ in }(0,T)\times \R^d\times \R^d,\\
    u(0,x,\vi)=g(x,\vi) \text{ in }\R^d\times \R^d,
    \end{cases}
\end{equation*}
according to Definition \ref{RenormHJBDefinition}.
\end{theorem}
\begin{proof}
Following \cite{porretta1999existence}, we see that the result will hold once we show that for some increasing sequence $0\leq m_k\in \R, k\in \N$ with $m_k\uparrow \infty$ as $k\rightarrow \infty$, $\Dv (T_{m_k}(\un))\rightarrow \Dv (T_{m_k}(u))$ strongly in $L^2([0,T]\times \R^d\times \R^d)$, where
\begin{equation}
    T_k(s):=\begin{cases}
    s, \text{ if }|s|\leq k,\\
    k, \text{ if }s>k,\\
    -k, \text{ if }s<-k.
    \end{cases}
\end{equation}
Note that for almost all $\beta\in \R$, we have that $|\{u=\beta\}|=0$ ($|A|$ denotes the Lebesgue measure), therefore in order to keep the notation lighter we may assume that $|\{u=k\}|=0$ and thus choose the sequence $m_k=k$. The reason for this choice will become apparent later; in particular to prove that $\chi_{u^n>m_k}\rightarrow \chi_{u>m_k}$ almost everywhere, it is convenient to know that $|\{u=m_k\}|=0$. In the rest of the proof we will use the notation $\omega(n)$ and $\omega(n,\alpha)$, for quantities that satisfy $\lim\limits_{n\rightarrow \infty}\omega(n)=0$ and $\lim\limits_{\alpha\rightarrow \infty}\lim\limits_{n\rightarrow\infty}\omega(n,\alpha)=0$ respectively, furthermore these quantities are subject to change from line to line. Just as in \cite{porretta1999existence} and the references therein, for $\lambda>0$ we define $\phil(s):=s\exp(\lambda s^2)$. For $\alpha >0$ and $k\in \N$, consider the solution $\Phiak$ of
\begin{equation}
    \begin{cases}
    \pt \Phiak-\lapv \Phiak+\vDx \Phiak=\alpha(T_k(u)-\Phiak),\\
    \Phiak(0)=T_k(g).
    \end{cases}
\end{equation}
Denote by $L:=\pt -\lapv +\vDx$ and test equation \eqref{HJBInPhiaConvergence} against $\phil(\un -\Phiak)^-$ which yields
\[\boxed{\int_0^T\intr \langle L(\un -\Phiak), \phil(\un -\Phiak)^-\rangle dxd\vi dt}_1\]
\[+\boxed{\int_0^T \intr \langle L\Phiak, \phil(\un -\Phiak)^-\rangle dxd\vi dt}_2\]
\[+\boxed{\int_0^T\intr H(\Dv \un )\phil(\un -\Phiak)^-dxd\vi dt}_3=\boxed{\int_0^T\intr f^n \phil(\un -\Phiak)^-dxd\vi dt}_4.\]
Let $\Phi_{\lambda}(s):=\int_0^s \phi_l(\theta)^-d\theta$, then the first boxed term gives us
\[\intr \Phi_{\lambda}(\un -\Phiak)(T) \Phi_{\lambda}(g^n-T_k(g))dxd\vi -\int_0^T\intr \phil'(\un -\Phiak)^-|\Dv (\un -\Phiak)|^2dxd\vi dt\]
\[\leq \omega(n)-\int_0^T\intr \phil'(\un -\Phiak)^-|\Dv (\un -\Phiak)|^2dxd\vi dt,\]
where in the last inequality we used that $\Phi_{\lambda}(s):=\int_0^s \phil(u)^-du\leq 0$ and that $g^n\rightarrow g$ strongly in $L^1$. For the second boxed term we obtain  
\[\alpha \int_0^T\intr (T_k(u)-\Phiak)\phil(\un -\Phiak)^-dxd\vi dt\leq \alpha\omega(n),\]
since $\un\rightarrow u$ strongly in $L^1$, $\phil(\un-\Phiak)^-=\phil(T_k(\un)-\Phiak)^-$ and $s\phil(s)^-\leq 0$. For the third boxed term we have that for some constant $C>0$, depending only on $H$
\[\int_0^T\intr H(\Dv \un )\phil(\un -\Phiak)^-dxd\vi dt\leq C\int_0^T\intr |\Dv (\un)|^2 \phil(\un -\Phiak)^-dxd\vi dt,\]
and using that for all $p,q\in \R^d$, we have $|p|^2\leq 2|p-q|^2+2|q|^2$ the third boxed term is bounded by  
\[2C\int_0^T\intr |\Dv (\un-\Phiak)|^2 \phil(\un -\Phiak)^- + 2C |\Dv (\Phiak)|^2 \phil(\un -\Phiak)^-dxd\vi dt.\]
Finally, the last boxed term vanishes as $n\rightarrow \infty$ and then $\alpha\rightarrow \infty$ due to Lemma \ref{PhiaConvergence}. Putting everything together we obtain 
\[2C\int_0^T\intr \Big[\phil'(\un -\Phiak)^--\phil(\un -\Phiak)^-\Big]|\Dv (\un -\Phiak)|^2dxd\vi dt\]
\[\leq \omega(n,\alpha)+2C\int_0^T\intr |\Dv (\Phiak)|^2 \phil(\un -\Phiak)^-dxd\vi dt.\]
By choosing $\lambda$ large enough depending only on $\|H_{pp}\|_{\infty},$ we have that $\phil'(\un -\Phiak)^--\phil(\un -\Phiak)^-\geq 0$ thus by Fatous Lemma on the LHS and the strong convergence of $\un \rightarrow u$ in $L^1$, as $n\rightarrow \infty$ we obtain
\[\int_0^T\intr \Big[\phil'(u -\Phiak)^--2C\phil(u -\Phiak)^-\Big]|\Dv (u -\Phiak)|^2dxd\vi dt\]
\[\leq \omega(\alpha)+2C\int_0^T\intr |\Dv (\Phiak)|^2 \phil(u -\Phiak)^-dxd\vi dt.\]
Furthermore,
\[2C\int_0^T\intr |\Dv (\Phiak)|^2 \phil(u -\Phiak)^-dxd\vi dt\]
\[\leq 4C\int_0^T\intr |\Dv (\Phiak-u)|^2 \phil(u -\Phiak)^-dxd\vi dt+4C\int_0^T\intr |\Dv u|^2 \phil(u -\Phiak)^-dxd\vi dt .\]
Hence,
\[\int_0^T\intr \Big[\phil'(u -\Phiak)^--6C\phil(u -\Phiak)^-\Big]|\Dv (u -\Phiak)|^2dxd\vi dt\]
\[\leq \omega(\alpha)+ 4C\int_0^T\intr |\Dv u|^2 \phil(u -\Phiak)^-dxd\vi dt,\]
now we may fix $\lambda>0$ such that $\phil'(s)^--6C\phil(s)^-\geq \frac{1}{2}$ and so letting $\alpha\rightarrow \infty$ yields  
\[\lim\limits_{\alpha\rightarrow \infty}\|\Dv (T_k(u)-\Phiak)^-\|_2=0.\]
We now show the convergence on the set $T_k(u)\geq \Phiak$. Since $H\geq 0$ the functions $\un $ are subsolutions of 
\begin{equation}\label{TruncationSubsolutionLemma}
    \begin{cases}
    L\un \leq f^n(t,x,\vi) \text{ in }(0,T)\times \R^d\times \R^d\\
    \un(0,x,\vi)=g^n(x,\vi).
    \end{cases}
\end{equation}
Define $w^n=(T_k(\un)-\Phiak)_+$ which may also be written as
\[w_n=(\un-\Phiak)_+-(\un-T_k(\un)).\]
Indeed if $\un\leq k$ then 
\[(\un-\Phiak)_+-(\un-T_k(\un))=(\un-\Phiak)_+=(T_k(\un)-\Phiak)_+,\]
while if $\un>k$ since $0\leq\Phiak\leq k$
\[(\un-\Phiak)_+-(\un-T_k(\un))=\un-\Phiak-\un+k=k-\Phiak=T_k(\un)-\Phiak=(T_k(\un)-\Phiak)_+.\]
Thus testing against $w_n$ in equation \eqref{TruncationSubsolutionLemma} yields
\[\int_0^T\intr \langle L(\un),w_n\rangle dxd\vi dt\leq \int_0^T\intr f^n w_ndxd\vi dt\implies \]
\[\int_0^T\intr\boxed{ \langle L(\un-\Phiak),(T_k(\un)-\Phiak)_+\rangle }_1+\boxed{\langle L(\Phiak),(T_k(\un)-\Phiak)_+\rangle dxd\vi dt}_2\]
\[\boxed{-\int_0^T\intr \langle L(\un),\un-T_k(\un)\rangle dxd\vi dt}\leq \boxed{\int_0^T\intr f^nw_ndxd\vi dt}_4.\]
The first boxed term equals
\[\int_0^T\intr \langle L(\un-\Phiak),(\un-\Phiak)_+\rangle dxd\vi dt=\intr (\un-\Phiak)_+^2/2(T)-((g^n-T_k(g))_+^2/2dxd\vi \]
\[+\int_0^T\intr \Dv(\un-\Phiak)\Dv (\un-\Phiak)_+dxd\vi dt\]
and since $g^n\in L^1\cap L^{\infty}$ the quantities that appear make sense. The second boxed term is bounded by 
\[\int_0^T\intr \langle L(\Phiak),(\un-\Phiak)_+\rangle dxd\vi dt=\alpha\int_0^T(T_k(u)-\Phiak)(\un-\Phiak)\geq \omega(n).\]
The third boxed term equals
\[-\int_0^T\intr \langle L(\un),\un-T_k(\un)\rangle dxd\vi dt\]
\[=-\int_0^T\intr (\un(T)-T_k(\un)(T))^2/2-(g^n-T_k(g^n))^2/2 \Dv \un\Dv(\un-T_k(\un))dxd\vi dt.\]
Putting everything together yields
\[\int_0^T\intr \langle L(\un),w_n\rangle dxd\vi dt\geq \omega(n)+\intr (\un-\Phiak)_+^2/2(T)-((g^n-T_k(g))_+^2/2dxd\vi\]
\[-\int_0^T\intr (\un(T)-T_k(\un)(T))^2/2-(g^n-T_k(g^n))^2/2dxd\vi\]
\[+\int_0^T\intr \Dv(\un-\Phiak)\Dv (\un-\Phiak)_+dxd\vi dt-\int_0^T\intr \Dv \un\Dv(\un-T_k(\un))dxd\vi dt .\]
The first line equals
\[\intr (\un-\Phiak)_+^2/2(T)-((g^n-T_k(g))_+^2/2- (\un(T)-T_k(\un)(T))^2/2-(g^n-T_k(g^n))^2/2dxd\vi\]
\[=\frac{1}{2}\intr \Big((\un-\Phiak)_+(T)-(\un(T)-T_k(\un)(T))((\un-\Phiak)_++(\un-T_k(\un))(T)) \Big)dxd\vi \]
\[-\frac{1}{2}\intr ((g^n-T_k(g))_+-(g^n-T_k(g^n)))((g^n-T_k(g))_++(g^n-T_k(g^n)))\]
\[\geq -2\frac{1}{2}\intr (T_k(g^n)-T_k(g))_+(g^n-T_k(g))_+dxd\vi =\omega(n).\]
For the last line
\[\int_0^T\intr \Dv(\un-\Phiak)\Dv (\un-\Phiak)_+dxd\vi dt-\int_0^T\intr \Dv \un\Dv(\un-T_k(\un))dxd\vi dt=\]
\[\int_0^T\intr \Dv(\un-\Phiak)\Dv\Big( (\un-\Phiak)_+-(\un-T_k(\un))\Big)+ \Dv(\un-\Phiak)\Dv (\un-T_k(\un))_+dxd\vi dt\]
\[-\int_0^T\intr \Dv \un\Dv(\un-T_k(\un))dxd\vi dt\]
\[=\int_0^T\intr \Dv(\un-\Phiak)\Dv(T_k(\un)-\Phiak)_+dxd\vi dt-\int_0^T\intr \Dv \Phiak \Dv(\un-T_k(\un))dxd\vi dt\]
\[=\int_0^T\intr |\Dv(\un-\Phiak)_+|^2+ \Dv(\un-T_k(\un))\Dv\Big((T_k(\un) -\Phiak)_+-\Phiak\Big)dxd\vi dt\]
\[=\int_0^T\intr |\Dv(T_k(\un)-\Phiak)_+|^2dxd\vi dt-2\int_0^T\int_{\un>k} \Dv(\un)\Dv\Big(\Phiak\Big)dxd\vi dt ,\]
where in the last equality we used that $\Dv(\un-T_k(\un))=\Dv \un\chi_{\un>k}$ and $0\leq \Phiak\leq k$. Finally, we clearly have that
\[\int_0^T\intr f^nw_ndxd\vi dt\leq \omega(n,\alpha).\]
Hence, putting everything together
\[\int_0^T\intr |\Dv(T_k(\un)-\Phiak)_+|^2dxd\vi dt\leq 2\int_0^T\int_{\un>k} \Dv(\un)\Dv\Big(\Phiak\Big)dxd\vi dt+\omega(n,\alpha).\]
Since $\Dv \un\rightarrow \Dv u$ weakly in $L^2$ while $\chi_{\un>k}\Phiak\rightarrow \chi_{u>k}\Phiak$ strongly in $L^2$ (here is where the discussion in the beginning of the proof is relevant) we may use Fatous Lemma which yields
\[\int_0^T\intr |\Dv(T_k(u)-\Phiak)_+|^2dxd\vi dt\leq 2\int_0^T\int_{u>k} \Dv(u)\Dv\Big(\Phiak\Big)dxd\vi dt+\omega(\alpha).\]
Furthermore,
\[\|\Dv \Phiak\|_2\leq \|\Dv(T_k(u)-\Phiak)_+\|_2+\|\Dv(T_k(u)-\Phiak)_-\|_2+\|\Dv T_k(u)\|_2\leq C,\]
for some $C>0$ independent of $\alpha$ (due to $\omega(\alpha)\rightarrow 0$ as $\alpha\rightarrow \infty$). Therefore, we may assume WLOG that $\Dv \Phia\rightarrow \Dv T_k(u)$ weakly in $L^2$. Thus, taking the limit as $\alpha\rightarrow\infty$ we find that 
\[\limsup\limits_{\alpha\rightarrow \infty}\int_0^T\intr |\Dv(T_k(u)-\Phiak)_+|^2dxd\vi dt\leq\]
\[\lim\limits_{\alpha\rightarrow\infty}\Big(2\int_0^T\int_{u>k} \Dv(u)\Dv\Big(\Phiak\Big)dxd\vi dt+\omega(\alpha)\Big)=2\int_0^T\int_{u>k} \Dv(u)\Dv\Big(\Dv T_k(u)\Big)dxd\vi dt=0.\]
Now that we have $\Dv \Phiak\rightarrow \Dv T_k(u)$ strongly in $L^2$, we may conclude since by the previous estimates
\[\int_0^T\intr |\Dv(T_k(u^n)-\Phiak)_+|^2dxd\vi dt\leq \omega(n,\alpha).\]

\end{proof}
We conclude this subsection with a sketch of the proof for the upper bound in Theorem \ref{SolutionToFPequation}. We recall the Fractional Gagliardo-Niremberg inequality.
\begin{proposition}\label{FractionalGagliardoNirembergProp}
(Fractional Gagliardo-Niremberg inequality). Let $z\in H^s(\R^d\times \R^d),$ where $s>0$. If $\theta\in (0,1)$ $p\in (1,\infty)$ are such that 
$$\theta\Big(\frac{1}{2}-\frac{s}{d}\Big)+\frac{1-\theta}{2}=\frac{1}{p}\iff \frac{1}{p}=\frac{1}{2}-\frac{\theta s}{d},$$
then
$$\|z\|_p\leq C\|D^sz\|_2^{\theta}\|z\|_2^{1-\theta},$$
where $D^s z_a=(\Dv^s z_a,D_x^s z_a)$.
\end{proposition}
\begin{corollary}\label{CorGagliardo}
Let $z\in L^2((0,T);H^s(\R^d\times \R^d))$. Then, for $p=2(1+\frac{2s}{d})$ and $\theta p=2$, we have
\[\Big(\int_0^T \|z(t)\|_p^pdt \Big)^{1/p}\leq \sup\limits_{t\in [0,T]}\|z(t)\|_2^{1-\theta}\|D^s z\|_2^{2/p}=\sup\limits_{t\in [0,T]}\|z(t)\|_2^{1-\theta}\|D^s z\|_{\Lt}^{\theta}.\]
\end{corollary}
\begin{proposition}\label{prop: UpperBoundForm}
Let $b\in L^{\infty}([0,T]\times \R^d\times \R^d)$ and $m_0$ a density which satisfies \ref{AssumptionInitial}. Furthermore, let $m\in C([0,T]; L^2(\R^d\times \R^d))$ be the distributional solution to 
\begin{equation}
    \begin{cases}
        \pt m -\lapv m +\vDx m -\divv(m b) =0 \text{ in }(0,T)\times \R^d\times \R^d,\\
        m(0,x,\vi)=m_0(x,\vi) \text{ in }\R^d\times \R^d.
    \end{cases}
\end{equation}
Then, there exists a $C_0=C_0(\|b\|_{\infty},T,\|m_0\|_2,\|m_0\|_{\infty})>0$, such that
\[\|m\|_{\infty}\leq C_0.\]
\end{proposition}
\begin{proof}
    The proof follows the work of F. Golse, C. Imbert, C. Mouhot and A. Vasseur in \cite{golse2016harnack}. To simplify the notation we define the operator 
    \[\mathcal{L}^*m := \pt m -\lapv m +\vDx m -\divv(m b).\]
    Furthermore, to reduce the technical steps we make the following reduction.
    By linearity it is enough to show the result in the case of $\|m_0\|_{\infty}\leq 1$. Moreover, we assume that $b$ is smooth with compact support, since the general case may be handled by approximation. We note that once $b$ is smooth and compactly supported, the density $m$ is bounded above, however this bound depends on $\|\divv(b)\|_{\infty}$. Nonetheless, at the level of smooth $b$ the functions $m,m^2$ are integrable. For $\alpha >1\geq \|m_0\|_{\infty}$ we set $m_{\alpha}:=(m-\alpha)_+$. Then, we have that $m_{\alpha}$ is a subsolution of 
    \begin{equation}\label{subsolm_a^2_1}
        \mathcal{L}^*m_{\alpha}-(1+\alpha)\mathbf{1}_{m>\alpha}\divv(b)\leq 0, m_{\alpha}(0)=0.
    \end{equation} 
    Moreover, for technical reasons we will also require the function $m_{\alpha}^2$, which is a subsolution of
\begin{equation}
\label{eq: subsolm_a^2_1}
\pt m_{\alpha}^2-\lapv m_{\alpha}^2-\vDx m_{\alpha}^2-\divv(m_{\alpha}^2b)-m_{\alpha}^2\divv(b)-2\alpha m_{\alpha}\divv(b)\leq 0, \,\,\, m_{\alpha}^2(0)=0,
\end{equation}
or equivalently
\begin{equation}
\label{subsolm_a^2_2}
\pt m_{\alpha}^2-\lapv m_{\alpha}^2-\vDx m_{\alpha}^2-2\divv(m_{\alpha}^2b)+\Dv(m_{\alpha}^2)\cdot b-2\alpha \divv(m_{\alpha} b)+2\alpha Dv(m_{\alpha})\cdot b\leq 0, m_{\alpha}^2(0)=0.
\end{equation}
The typical energy estimates required in the De Giorgi argument for improvement of integrability, are not suitable for this setting. Namely testing against $m_a^2$ in \ref{eq: subsolm_a^2_1}, only yields bounds on $\Dv m_a^2$. To obtain an increase in integrability we first look at the solution $w_{\alpha}$ of
\begin{equation}
\label{equationw_a10}
\pt w_{\alpha}-\lapv w_{\alpha}-\vDx w_{\alpha}-\divv(m_{\alpha}^2b)-m_{\alpha}^2\divv(b)-2\alpha m_{\alpha}\divv(b)=0, \,\,\, w_{\alpha}(0)=0, 
\end{equation}
and we note that $w_a\geq m_a^2\geq 0$. Testing against $w_a$ in \eqref{equationw_a10} yields by Gr\"onwall
\begin{equation}
\begin{aligned}
\label{estimate_w_a_12}
&\sup\limits_{t\in [0,T]}\|w_{\alpha}(t)\|_2^2+\|\Dv w_{\alpha}\|_{L^2([0,T]\times \R^d\times \R^d)}^2\\
\leq C(\|m_{\alpha}\|_{\Lt}^2+&\|m_{\alpha}^2\|_{\Lt}^2+\|\Dv m_{\alpha}\|_{\Lt}^2+\|\Dv (m_{\alpha}^2)\|_{\Lt}^2).
\end{aligned}
\end{equation}
For the estimates on $m_{\alpha}^2$ we test \eqref{subsolm_a^2_1} against $m_{\alpha}^2$ and integrate in space to obtain by Gr\"onwall
\begin{equation}
\label{Ineq13}
\sup\limits_{t\in [0,T]}\int |m_{\alpha}(t)|^4+\int_0^T \int |\Dv m_{\alpha}^2|^2\leq C(\int_0^T \int |m_{\alpha}|^2+\int_0^T\int |\Dv m_{\alpha}|^2).
\end{equation}
We need an estimate for $\int_0^T\int |\Dv m_{\alpha}|^2$, so we test against $m_{\alpha}$ in \eqref{subsolm_a^2_1} and integrate in space to obtain
by Gr\"onwall
\begin{equation}
\label{estimate14}
\sup\limits_{t\in [0,T]}\|m_{\alpha}(t)\|_2^2+\int_0^T\int |\Dv m_{\alpha}|^2\leq C\int_0^T |\{m_{\alpha}(t)>0\}|.
\end{equation}
Using estimates \eqref{estimate14},\eqref{Ineq13} on \eqref{estimate_w_a_12} yields
\begin{equation}
\label{estimate15}
\sup\limits_{t\in [0,T]}\|w_{\alpha}(t)\|_2^2+\int_0^{T+2}\int |\Dv w_{\alpha}|^2\leq C\int_0^T |\{m_{\alpha}(t)>0\}|dt.
\end{equation}
From the above and Theorem \ref{Thm2.1FiniteTime}, we obtain
\begin{equation}
\label{estimate16}
\|D^s w_{\alpha}\|_{\Lt}^2\leq C\int_0^T |\{m_{\alpha}(t)>0\}|dt.
\end{equation}
From \eqref{estimate16} and Corollary \ref{CorGagliardo}, we obtain
\[\|w_{\alpha}\|_{L^p([0,T]\times \R^d\times \R^d)}\leq C\|D^s w_{\alpha}\|_{\Lt}^{\theta}\sup\limits_{t\in [0,T]}\|w_{\alpha}(t)\|_{\Lt}^{1-\theta}\leq C\|D^s w_{\alpha}\|_2^{\theta}\sup\limits_{t\in [0,T]}\|w_{\alpha}(t)\|_2^{1-\theta}\]
from \eqref{estimate15} and \eqref{estimate16} we have
\begin{equation}
\label{estimate17}
\|w_{\alpha}\|_{L^p([0,T]\times \R^d\times \R^d)}\leq C\int_0^T |\{m_{\alpha}(t)>0\}|.
\end{equation}
We may now setup the De-Giorgi iteration. For $k\in \N$, let $\alpha_k=(2+\frac{1}{2^{k-1}})$ and $m_k:=m_{\alpha_k}$. Since
\begin{equation}
\label{estimate18}
|\{m_k(t)>0\}|=|\{m_{k-1}(t)>\frac{1}{2^k}\}|\leq 16^k\int |m_{k-1}(t)|^4,
\end{equation}
if we define $U_k:=\int_0^T \intr |m_k|^4dxd\vi dt$, and use estimate \eqref{estimate18} in \eqref{estimate17}, we obtain
\begin{equation}
\label{estimate19}
\|w_k\|_{L^p([0,T]\times \R^d\times \R^d)}\leq C16^kU_{k-1}.
\end{equation}
Recall that $m_{\alpha}^2\leq w_{\alpha},$ thus from \eqref{estimate19} we have
\[\|m_{\alpha}^2\|_p\leq \|w_{\alpha}\|_p\leq C16^kU_{k-1}.\]
Therefore, 
\[U_k=\int_0^T\int |m_k|^4dxd\vi dt=\|m_k^2\|_2^2\leq C\|w_k\|_p^2|\{m_k>0\}|^{\epsilon}\leq C16^kU_{k-1}^{1+\epsilon},\]
for some $\epsilon=\epsilon(p)>0$ and the result follows.
\end{proof}
\subsection{Prerequisites}
We rely on the following minor modifications of three results from \cite{bouchut2002hypoelliptic}. We modify these Theorems slightly, to be used for a finite time interval $[0,T]$. 
\begin{theorem}\label{Thm1.2FiniteTime}
(Theorem 1.5,\cite{bouchut2002hypoelliptic}) Let $f,g\in L^2([0,T]\times \R^d\times \R^d)$, $\Dv f\in L^2(\R\times \R^d\times \R^d)$ and $f_0\in L^2(\R^d\times \R^d)$, such that
\begin{equation*}
    \begin{cases}
    \pt f-\vDx f-\lapv f=g \text{ in }[0,T]\times \R^d\times \R^d,\\
    f(0,x,\vi)=f_0(x,\vi) \text{ in }\R^d\times \R^d.
    \end{cases}
\end{equation*}
Then, there exists a dimensional constant $C>0$, such that
\[\|\pt f-\vDx f\|_2+\|\lapv f\|_2\leq \frac{C}{t}\Big(\|g\|_2+\|f_0\|_2\Big).\]
\end{theorem}
\begin{theorem}\label{Thm1.3FiniteTime}
(Theorem 1.3, \cite{bouchut2002hypoelliptic}) Assume that $f,g,g_0\in L^p([0,T]\times \R^d\times \R^d)$, with $\Dv f\in L^p([0,T]\times \R^d\times \R^d)$, $(1+|\vi|^2)g\in L^p(\R\times \R^d\times \R^d)$, $(1+|\vi|)g_0\in L^p(\R\times \R^d\times \R^d)$ and $f_0\in L^p(\R^d\times \R^d)$ for some $p\in (1,\infty)$, such that they solve
\begin{equation*}
\begin{cases}
\pt f-\vDx f=\divv(g)+g_0 \text{ in }(0,T]\times \R^d\times \R^d,\\
f(0,x,\vi)=f_0(x,\vi) \text{ in }\R^d\times \R^d,
\end{cases}
\end{equation*}
in the distributional sense. Then, there exists a constant $C>0$, such that
\[\|D_x^{1/3}f\|_p+\|D_t^{1/3}f\|_p\leq C(\|f\|_p+\|\Dv f\|_p+\|(1+|\vi|^2)g\|_p+\|(1+|\vi|)g_0\|_p+\|f_0\|_p).\]
\end{theorem}
\begin{theorem}\label{Thm2.1FiniteTime}
(Theorem 2.1, \cite{bouchut2002hypoelliptic}) Assume that $f,g,g_0\in L^p([0,T]\times \R^d\times \R^d)$, with $\Dv f\in L^p([0,T]\times \R^d\times \R^d)$, and $f_0\in L^p(\R^d\times \R^d)$ for some $p\in (1,\infty)$, such that they solve
\begin{equation*}
\begin{cases}
\pt f-\vDx f=\divv(g)+g_0 \text{ in }(0,T]\times \R^d\times \R^d,\\
f(0,x,\vi)=f_0(x,\vi) \text{ in }\R^d\times \R^d,
\end{cases}
\end{equation*}
in the distributional sense . Then, there exists a constant $C>0$, such that
\[\|D_x^{1/3}f\|_p\leq C(\|\Dv f\|_p+\|f_0\|_p+\|g\|_p+\|g_0\|_p),\]
where $\alpha,\alpha'\in (0,1)$ and depend only on the dimension d.
\end{theorem}
\section*{Acknowledgement}
The author was  partially supported by P. E. Souganidis' NSF grant DMS-1900599, ONR grant N000141712095 and AFOSR  grant FA9550-18-1-0494. 
\bibliographystyle{abbrv}
\bibliography{references}

\begin{thebibliography}{10}

\bibitem{achdou2020deterministic}
Y.~Achdou, P.~Mannucci, C.~Marchi, and N.~Tchou.
\newblock Deterministic mean field games with control on the acceleration.
\newblock {\em Nonlinear Differential Equations and Applications NoDEA},
  27(3):1--32, 2020.

\bibitem{bardi2021convergence}
M.~Bardi and P.~Cardaliaguet.
\newblock Convergence of some mean field games systems to aggregation and
  flocking models.
\newblock {\em Nonlinear Analysis}, 204:112199, 2021.

\bibitem{boccardo1997nonlinear}
L.~Boccardo, A.~Dall'Aglio, T.~Gallou{\"e}t, and L.~Orsina.
\newblock Nonlinear parabolic equations with measure data.
\newblock {\em journal of functional analysis}, 147(1):237--258, 1997.

\bibitem{bouchut2002hypoelliptic}
F.~Bouchut.
\newblock Hypoelliptic regularity in kinetic equations.
\newblock {\em Journal de math{\'e}matiques pures et appliqu{\'e}es},
  81(11):1135--1159, 2002.

\bibitem{brezis}
H.~Brezis and H.~Br{\'e}zis.
\newblock {\em Functional analysis, Sobolev spaces and partial differential
  equations}, volume~2.
\newblock Springer, 2011.

\bibitem{brezis2018gagliardo}
H.~Brezis and P.~Mironescu.
\newblock Gagliardo--nirenberg inequalities and non-inequalities: the full
  story.
\newblock In {\em Annales de l'Institut Henri Poincar{\'e} C, Analyse non
  lin{\'e}aire}, volume~35, pages 1355--1376. Elsevier, 2018.

\bibitem{burger2016balanced}
M.~Burger, A.~Lorz, and M.-T. Wolfram.
\newblock Balanced growth path solutions of a boltzmann mean field game model
  for knowledge growth.
\newblock {\em arXiv preprint arXiv:1602.01423}, 2016.

\bibitem{camilli2021quadratic}
F.~Camilli.
\newblock A quadratic mean field games model for the langevin equation.
\newblock {\em Axioms}, 10(2):68, 2021.

\bibitem{cardaliaguet2010notes}
P.~Cardaliaguet.
\newblock Notes on mean field games.
\newblock Technical report, Technical report, 2010.

\bibitem{cardaliaguet2015second}
P.~Cardaliaguet, P.~J. Graber, A.~Porretta, and D.~Tonon.
\newblock Second order mean field games with degenerate diffusion and local
  coupling.
\newblock {\em Nonlinear Differential Equations and Applications NoDEA},
  22(5):1287--1317, 2015.

\bibitem{carmona2018probabilistic}
R.~Carmona and F.~Delarue.
\newblock Probabilistic theory of mean field games with applications. i, volume
  83 of probability theory and stochastic modelling, 2018.

\bibitem{degond1986global}
P.~Degond.
\newblock Global existence of smooth solutions for the vlasov-fokker-planck
  equation in $1 $ and $2 $ space dimensions.
\newblock In {\em Annales scientifiques de l'{\'E}cole Normale Sup{\'e}rieure},
  volume~19, pages 519--542, 1986.

\bibitem{demengel2012fractional}
F.~Demengel, G.~Demengel, F.~Demengel, and G.~Demengel.
\newblock Fractional sobolev spaces.
\newblock {\em Functional spaces for the theory of elliptic partial
  differential equations}, pages 179--228, 2012.

\bibitem{diperna1988fokker}
R.~J. DiPerna and P.-L. Lions.
\newblock On the fokker-planck-boltzmann equation.
\newblock {\em Communications in mathematical physics}, 120(1):1--23, 1988.

\bibitem{dragoni2018ergodic}
F.~Dragoni and E.~Feleqi.
\newblock Ergodic mean field games with h{\"o}rmander diffusions.
\newblock {\em Calculus of Variations and Partial Differential Equations},
  57(5):1--22, 2018.

\bibitem{feleqi2020hypoelliptic}
E.~Feleqi, D.~A. Gomes, and T.~Tada.
\newblock Hypoelliptic mean-field games—a case study.
\newblock 2020.

\bibitem{folland1975subelliptic}
G.~B. Folland.
\newblock Subelliptic estimates and function spaces on nilpotent lie groups.
\newblock {\em Arkiv f{\"o}r matematik}, 13(1):161--207, 1975.

\bibitem{golse2016harnack}
F.~Golse, C.~Imbert, C.~Mouhot, and A.~Vasseur.
\newblock Harnack inequality for kinetic fokker-planck equations with rough
  coefficients and application to the landau equation.
\newblock {\em arXiv preprint arXiv:1607.08068}, 2016.

\bibitem{griffin2022variational}
M.~Griffin-Pickering and A.~R. M{\'e}sz{\'a}ros.
\newblock A variational approach to first order kinetic mean field games with
  local couplings.
\newblock {\em Communications in Partial Differential Equations},
  47(10):1945--2022, 2022.

\bibitem{huang2006large}
M.~Huang, R.~P. Malham{\'e}, P.~E. Caines, et~al.
\newblock Large population stochastic dynamic games: closed-loop mckean-vlasov
  systems and the nash certainty equivalence principle.
\newblock {\em Communications in Information \& Systems}, 6(3):221--252, 2006.

\bibitem{lasry2006jeux}
J.-M. Lasry and P.-L. Lions.
\newblock Jeux {\`a} champ moyen. i--le cas stationnaire.
\newblock {\em Comptes Rendus Math{\'e}matique}, 343(9):619--625, 2006.

\bibitem{lasry2006_2jeux}
J.-M. Lasry and P.-L. Lions.
\newblock Jeux {\`a} champ moyen. ii--horizon fini et contr{\^o}le optimal.
\newblock {\em Comptes Rendus Math{\'e}matique}, 343(10):679--684, 2006.

\bibitem{lasry2007mean}
J.-M. Lasry and P.-L. Lions.
\newblock Mean field games.
\newblock {\em Japanese journal of mathematics}, 2(1):229--260, 2007.

\bibitem{lunardi1997schauder}
A.~Lunardi.
\newblock Schauder estimates for a class of degenerate elliptic and parabolic
  operators with unbounded coefficients in $\mathbb{R}^n$.
\newblock {\em Annali della Scuola Normale Superiore di Pisa-Classe di
  Scienze}, 24(1):133--164, 1997.

\bibitem{mouhot2018giorgi}
C.~Mouhot.
\newblock De giorgi--nash--moser and h{\"o}rmander theories: new interplays.
\newblock In {\em Proceedings of the International Congress of Mathematicians
  (ICM 2018) (In 4 Volumes) Proceedings of the International Congress of
  Mathematicians 2018}, pages 2467--2493. World Scientific, 2018.

\bibitem{porretta1999existence}
A.~Porretta.
\newblock Existence results for nonlinear parabolic equations via strong
  convergence of truncations.
\newblock {\em Annali di Matematica Pura ed Applicata}, 177(1):143--172, 1999.

\bibitem{porretta2015weak}
A.~Porretta.
\newblock Weak solutions to fokker--planck equations and mean field games.
\newblock {\em Archive for Rational Mechanics and Analysis}, 216(1):1--62,
  2015.

\bibitem{porretta2017weak}
A.~Porretta.
\newblock On the weak theory for mean field games systems.
\newblock {\em Bollettino dell'Unione Matematica Italiana}, 10(3):411--439,
  2017.

\bibitem{Adjoint}
E.~Priola.
\newblock Global schauder estimates for a class of degenerate kolmogorov
  equations.
\newblock {\em arXiv preprint arXiv:0705.2810}, 2007.

\bibitem{simon1986compact}
J.~Simon.
\newblock Compact sets in the space $l^p((0,t);b)$.
\newblock {\em Annali di Matematica pura ed applicata}, 146(1):65--96, 1986.

\end{thebibliography}

\end{document}